\documentclass{article}
\usepackage[T2A]{fontenc}
\usepackage[cp1251]{inputenc}
\usepackage[english,russian]{babel}
\usepackage[tbtags]{amsmath}
\usepackage{amsfonts,amssymb,mathrsfs,amscd}
\numberwithin{equation}{section}
\newtheorem{remark}{Remark}[section]
\newtheorem{lemma}{Lemma}[section]
\newtheorem{theorem}{Theorem}[section]

\newtheorem{corollary}{Corollary}[section]

\renewcommand{\Re}{\mathop{\rm Re}}

\newcommand{\supp}{\mathop{\rm supp}}

\newcommand{\col}{\mathop{\rm col}}

\newcommand{\diag}{\mathop{\rm diag}}

\begin{document}
\begin{Large}
\thispagestyle{empty}
\begin{center}
{\bf Inverse scattering problem for a third-order operator with local potential\\
\vspace{5mm}
V. A. Zolotarev}\\

B. Verkin Institute for Low Temperature Physics and Engineering
of the National Academy of Sciences of Ukraine\\
47 Nauky Ave., Kharkiv, 61103, Ukraine

Department of Higher Mathematics and Informatics, V. N. Karazin Kharkov National University \\
4 Svobody Sq, Kharkov, 61077,  Ukraine

\end{center}
\vspace{5mm}

{\small {\bf Abstract.} Inverse scattering problem for the operator representing sum of the operator of the third derivative on semi-axis and of the operator of multiplication by a real function is studied in this paper. Properties of Jost solutions of such an operator are studied and it is shown that these Jost solutions are solutions of the Riemann boundary value problem on a system of rays. The main system of linear singular integral equations is derived. This system is equivalent to the solution of inverse scattering problem.}
\vspace{5mm}


\vspace{5mm}

\begin{center}
{\bf Introduction}
\end{center}
\vspace{5mm}

Inverse scattering problem for the Schr$\ddot{\rm o}$dinger operator is one of intriguing problems of mathematical physics. Method of solution of the inverse problem (in the case of one variable) is based on classic works by I. M. Gelfand -- B. M. Levitan, V. A. Marchenko, M. G. Krein \cite{1} -- \cite{5}. It is this method that became an important instrument for integration of non-linear differential equations \cite{1} -- \cite{6}.

Solution of the inverse problem for differential equations of order $p$ ($p>2$) met with significant difficulties which is partially explained by the absence of mapping operators when $p>2$. One ought to note that in the process of integration of non-linear equations describing wave propagation in dispersive media, viz., Camassa -- Holm and Degasperis -- Procesi equations \cite{7} -- \cite{12}, Lax pairs L - A arise in which operator $L$ is of third order (qubic string). Therefore there exists a natural need in study of the inverse problem for differential operators of third order. Inverse problems (spectral and scattering) for operators of third order with non-local potential are solved in \cite{13} -- \cite{15}. Note that in \cite{14, 15} the function classes (analogues of cosines and sines) are defined. These functions serve as a natural language (tool) in the study of direct and inverse problems.

This paper solves the inverse scattering problem for a third-order operator with local potential
$$(Ly)(x)=iD^3y(x)+q(x)y(x)\quad(y\in L^2(\mathbb{R}_+))$$
on the semi-axis ($x\in\mathbb{R}_+$) where $D=d/dx$ and $q(x)$ is a real function, such that
$$\int\limits_0^\infty|q(x)|^2e^{2ax}dx<\infty\quad(a>0).$$
We can take the `coefficient' of a symmetrized operator of first derivative (which replaces the term $q(x)y(x)$) as a potential. Operator $L$, domain of which consists of trice differentiable functions such that $y(0)=0$, $y'(0)=0$, is symmetric and does not have self-adjoint extensions into $L^2(\mathbb{R}_+)$ \cite{16} -- \cite{19}, and these extensions are constructed with exit from the space $L^2(\mathbb{R}_+$). For scattering problem, one has to extend $L$ up to a self-adjoint operator \cite{19}, i. e., to add scattering channels providing unitarity of scattering matrix. Therefore properties of the scattering matrix depend on selected self-adjoint extension. System of linear singular equations, solvability of which is equivalent to the inverse scattering problem, is derived in this paper. This equation system is obtained from the analysis of the Riemann boundary value problem on a system of rays.

The paper consists of 4 sections. The first section states the main properties \cite{14, 15} of the fundamental system of solutions of the operator $iD^3$ (an analogue of cosines and sines). Fourier transforms on a system of rays are defined and classes of functions arising during this process are described.

The second section studies Jost solutions of the self-adjoint extension of the operator $L$. It is established that Jost solutions exist and are analytic (relative to the spectral parameter) in corresponding sectors. It is shown that zeros of Jost functions are located on two rays that divide sector of holomorphicity into three equal sectors.

The third section studies scattering problem. In this case, we have two coefficients of scattering since one incident wave generates two scattered waves. Study of properties of the coefficients leads to the Riemann boundary value problem on the system of three rays \cite{15, 20}. A closed system of singular integral equations in terms of the data of inverse problem is derived from this Riemann problem. Note that condition for unitarity of the scattering problem (matrix), that follows from the self-adjoint boundary conditions of extensions of the operator $L$, has a non-standard form.

The fourth section of the paper solves the inverse problem, i. e., it analyzes solvability of the main system of singular equations. Scattering data for a third-order operator is described.

\section{Preliminary information}\label{s1}

{\bf 1.1} Equation
\begin{equation}
iy'''(x)=\lambda^3y(x)\quad(x\in\mathbb{R},\,\lambda\in\mathbb{C})\label{eq1.1}
\end{equation}
has \cite{14, 15, 21} three linear independent solutions $\{e^{i\lambda\zeta_kx}\}_1^3$ where $\{\zeta_k\}_1^3$ are roots of the equation $\zeta^3=1$,
\begin{equation}
\zeta_1=1;\quad\zeta_2=\frac12(-1+i\sqrt3);\quad\zeta_3=\frac12(-1-i\sqrt3)\label{eq1.2}
\end{equation}
and any solution to equation \eqref{eq1.1} is a linear combination of $\{e^{i\lambda\zeta_kx}\}_1^3$. Instead of the exponents $\{e^{i\lambda\zeta_kx}\}$, it is convenient to take another system of fundamental solutions to the equation \eqref{eq1.1}, $\{s_p(i\lambda x)\}_0^2$, where
\begin{equation}
s_p(z)\stackrel{\rm def}{=}\frac13\sum\limits_{k=1}^3\frac1{\zeta_k^p}e^{z\zeta_k}\quad(p=0,1,2;\,z\in\mathbb{C}).\label{eq1.3}
\end{equation}
Functions $\{s_p(i\lambda x)\}_0^2$ are analogous to cosines and sines for the second-order equation ($-y''(x)=\lambda^2y(x)$).

\begin{lemma}{[15]}\label{l1.1}
Entire functions $\{s_p(z)\}_0^2$ \eqref{eq1.3} of exponential type satisfy the relations

{\rm (i)} $s'_2(z)=s_1(z);$ $s'_1(z)=s_0(z);$ $s'_0(z)=s_2(z);$

{\rm (ii)} reality, $\overline{s_p(z)}=s_p(\overline{z})$ $(0\leq p\leq 2);$

{\rm (iii)} $p$-evenness, $s_p(\zeta_2z)=z_2^ps_p(z)\quad(0\leq p\leq2);$

{\rm (iv)} the Euler formula
$$e^{z\zeta_k}=s_0(z)+\zeta_ks_1(z)+\zeta_k^2s_2(z)\quad(1\leq k\leq3);$$

{\rm (v)} functions $\{s_p(z)\}_0^2$, solutions to the equation $y'''(z)=y(z)$, satisfying the initial data
$$s_0(0)=1;\quad s'_0(0)=0;\quad s''_0(0)=0;$$
$$s_1(0)=0;\quad s'_1(0)=1;\quad s''_1(0)=0;$$
$$s_2(0)=0;\quad s'_2(0)=0;\quad s''_2(0)=1;$$

{\rm (vi)} the main identity
$$s_0^3(z)+s_1^3(z)+s_2^3(z)-3s_0(z)s_1(z)s_2(z)=1;$$

{\rm (vii)} addition formulas
$$s_0(z+w)=s_0(z)s_0(w)+s_1(z)s_2(w)+s_2(z)s_1(w);$$
$$s_1(z+w)=s_0(z)s_1(w)+s_1(z)s_0(w)+s_2(z)s_2(w);$$
$$s_2(z+w)=s_0(z)s_2(w)+s_1(z)s_1(w)+s_2(z)s_0(w);$$

{\rm (viii)} $3s_0(z)s_0(w)=s_0(z+w)+s_0(z+\zeta_2w)+s_0(z+\zeta_3w);$
$$3s_0(z)s_2(w)=s_0(z+w)+\zeta_2s_0(z+\zeta_2w)+\zeta_3s_0(z+\zeta_3w);$$
\begin{equation}
3s_0(z)s_1(w)=s_1(z+w)+s_1(z+\zeta_2w)+s_1(z+\zeta_3w);
\label{eq1.4}
\end{equation}
$$3s_2(z)s_2(w)=s_1(z+w)+\zeta_2s_1(z+\zeta_2w)+\zeta_3s_1(z+\zeta_3w);$$
$$3s_0(z)s_2(w)=s_2(z+w)+s_2(z+\zeta_2w)+s_2(z+\zeta_3w);$$
$$3s_1(z)s_1(w)=s_2(z+w)+\zeta_2s_2(z+\zeta_2w)+\zeta_3s_3(z+\zeta_3w);$$

{\rm (ix)} $3s_0^2(z)=s_0(2z)+2s_0(-z);$

$\quad\hspace{3mm} 3s_1^2(z)=s_2(2z)+2s_2(-z);$

$\quad\hspace{3mm} 3s_2^2(z)=s_1(2z)+2s_1(-z);$

{\rm (x)} $s_0^2(z)-s_1(z)s_2(z)=s_0(-z);$

$\quad \hspace{3mm} s_1(z)-s_2(z)s_0(z)=s_2(-z);$

$\quad \hspace{3mm} s_2^2(z)-s_0(z)s_1(z)=s_1(-z);$

{\rm (xi)} the Taylor formula

$\quad\hspace{3mm}{\displaystyle s_0(z)=1+\frac{z^3}{3!}+\frac{z^6}{6!}+...};$

$\quad\hspace{3mm}{\displaystyle s_1(z)=z+\frac{z^4}{4!}+\frac{z^7}{7!}+...};$

$\quad\hspace{3mm}{\displaystyle s_2(z)=\frac{z^2}{2!}+\frac{z^5}{5!}+\frac{z^8}{8!}+....}$
\end{lemma}

\begin{remark}\label{r1.1}
For the equation $(i)^{-n}y^{(n)}(x)=\lambda^ny(x)$, for an arbitrary $n\in\mathbb{N}$, analogously to \eqref{eq1.3}, it is easy to construct the fundamental system of equations $\{s_p(i\lambda x)\}_0^{n-1}$ for which analogue of Lemma \ref{l1.1} (in proper formulation) is true.
\end{remark}

Solution to the Cauchy problem
\begin{equation}
iy'''(x)=\lambda^3y(x)+f(x);\quad y(0)=y_0,\quad y'(0)=y_1,\quad y''(0)=y_2,\label{eq1.5}
\end{equation}
for $f=0$, due to (v) \eqref{eq1.4}, is
\begin{equation}
y_0(\lambda,x)\stackrel{\rm def}{=}y_0s_0(i\lambda x)+y_1\frac{s_1(i\lambda x)}{i\lambda}+y_2\frac{s_2(i\lambda x)}{(i\lambda)^2}.\label{eq1.6}
\end{equation}
Using method of variation of constants, we find the solution to the non-homogenous ($f\not=0$) Cauchy problem
\begin{equation}
y(\lambda,x)=y_0(\lambda,x)-i\int\limits_0^x\frac{s_2(i\lambda(x-t))}{(i\lambda)^2}f(t)dt\label{eq1.7}
\end{equation}
where $y_0(\lambda,x)$ is given by \eqref{eq1.6}.

Unit vectors $\{\zeta_k\}_1^3$ define three straight lines,
\begin{equation}
L_{\zeta_k}\stackrel{\rm def}{=}\{x\zeta_k:x\in\mathbb{R}\}\quad(1\leq k\leq3),\label{eq1.8}
\end{equation}
besides, $\zeta_2L_{\zeta_1}=L_{\zeta_2}$, $\zeta_2L_{\zeta_2}=L_{\zeta_3}$, $\zeta_2L_{\zeta_3}=L_{\zeta_1}$. By $l_{\zeta_k}$, denote the straight half-lines (rays from the origin in the direction of unit vectors $\zeta_k$) and by $\widehat{l}_{\zeta_k}$, their supplements on $L_{\zeta_k}$ (rays coming to the origin),
\begin{equation}
l_{\zeta_k}\stackrel{\rm def}{=}\{x\zeta_k:x\in\mathbb{R}_+\};\quad\widehat{l}_{\zeta_k}\stackrel{\rm def}{=}L_{\zeta_k}\backslash l_{\zeta_k}=\{x\zeta_k:x\in\mathbb{R}_-\}\quad(1\leq k\leq3).\label{eq1.9}
\end{equation}
The straight lines $\{L_{\zeta_k}\}_1^3$ \eqref{eq1.8} divide plane $\mathbb{C}$ into the six sectors:
\begin{equation}
S_p\stackrel{\rm def}{=}\left\{z\in\mathbb{C}:\frac{2\pi}6(p-1)<\arg z<\frac{2\pi}6p\right\}\quad(1\leq p\leq6).\label{eq1.10}
\end{equation}

\begin{lemma} {[14, 15]}\label{l1.2}
Zeros of the functions $\{s_p(z)\}_0^2$ lie on the rays $\{\widehat{l}_{\zeta_k}\}_1^3$ \eqref{eq1.9} and are given by $\{-\zeta_2^lx_p(k)\}_{k=1}^\infty$ where $l=-1$, $0$, $1$, and $x_p(k)$ are non-negative numbers enumerated in the ascending order.

Besides, $x_0(k)$ are positive simple roots of the equation
$$\cos\frac{\sqrt3}2x=-\frac12 e^{-3/2 x}\quad(x_0(k)>0)$$
and $\{x_1(x)\}$ and $\{x_2(x)\}$ are non-negative simple roots of the equations
$$\cos\left(\frac{\sqrt3}2 x-\frac\pi 3\right)=\frac12e^{-\frac32x}\,(x_1(1)=0);\quad\cos\left(\frac{\sqrt3}2x+\frac\pi 3\right)=\frac12e^{-\frac32x}\,(x_2(1)=0).$$

The sequence $\{x_2(k)\}$ interlaces with the sequence $\{x_1(k)\}$ which, in its turn, interlaces with the sequence $\{x_0(k)\}$.
\end{lemma}

The equations for zeros $\{x_p(k)\}$ ($p=0$, 1, 2) readily yield asymptotic behavior of $x_p(k)$ as $k\rightarrow\infty$.
\vspace{5mm}

{\bf 1.2} The Fourier transform \cite{22, 23}
\begin{equation}
\widetilde{f}(\lambda)=\mathcal{F}(f)(\lambda);\quad\widetilde{f}(\lambda)\stackrel{\rm def}{=}\int\limits_{\mathbb{R}}e^{-i\lambda x}f(x)dx\quad(f\in L^2(\mathbb{R})),\label{eq1.11}
\end{equation}
due to $f=f_-+f_+$ ($f_\pm=f\chi_\pm$, $\chi_\pm$ are characteristic functions of the semiaxes $\mathbb{R}_\pm$) is given by
$$\widetilde{f}(\lambda)=\int\limits_0^\infty e^{-i\lambda x}f_+(x)dx+\int\limits_0^\infty e^{i\lambda x}f_-(-x)dx.$$
So, $\widetilde{f}(\lambda)$ is equal to the sum of Fourier transforms on the rays $\mathbb{R}_\pm$,
\begin{equation}
\widetilde{f}(\lambda)=\widetilde{f}_+(\lambda)+\widetilde{f}_-(\lambda);\quad\widetilde{f}_\pm(\lambda)\stackrel{\rm def}{=}\int\limits_0^\infty e^{-i\lambda\zeta_\pm x}f_\pm(\zeta_\pm x)dx\label{eq1.12}
\end{equation}
where $\zeta_\pm=\pm 1$ and $\{e^{-i\lambda\zeta_\pm x}\}$ is the system of linearly independent solutions to the equation $-D^2y(x)=\lambda^2y(x).$ Functions $\widetilde{f}_\pm(\lambda)$ \eqref{eq1.12} belong to $L^2(\mathbb{R})$ and $\widetilde{f}_+\perp\widetilde{f}_-$, besides, Parseval's identity $\|f_\pm\|_{L^2(\mathbb{R}_\pm)}=2\pi\|\widetilde f_\pm\|_{L^2(\mathbb{R})}$ holds \cite{22, 23}. Due to Paley -- Wiener theorem \cite{22, 23}, the function $\widetilde{f}_+(\lambda)$ ($\widetilde{f}_-(\lambda)$) is holomorphically extendable into the lower $\mathbb{C}_-$ (upper $\mathbb{C_+}$) half-plane and is of Hardy class $H_-^2$ ($H_+^2$) corresponding to $\mathbb{C}_-$ ($\mathbb{C}_+$). Thus, $\widetilde{f}(\lambda)$ \eqref{eq1.11} can be naturally represented by Fourier transform on the rays \eqref{eq1.12} where $\{e^{-i\lambda\zeta_\pm x}\}$ is the fundamental system of solutions to the second-order equation.

Describe the Fourier transform generated by the system of three exponents $\{e^{i\lambda\zeta_kx}\}_1^3$ where $\{\zeta_k\}$ are given by \eqref{eq1.2}. This system forms the fundamental system of solutions to the third-order equation \eqref{eq1.1}. Consider the bundle
\begin{equation}
l\stackrel{\rm def}{=}\bigcup\limits_kl_{\zeta_k}\label{eq1.13}
\end{equation}
formed by the rays $l_{\zeta_k}$ \eqref{eq1.9}, and let $\chi_k=\chi_{l_{\zeta_k}}$ be characteristic functions of the sets $l_{\zeta_k}$ ($1\leq k\leq 3$). Denote by $L^2(l)$ a Hilbert function space,
\begin{equation}
\begin{array}{ccc}
L^2(l)\stackrel{\rm def}{=}\{f=f_1\chi_1+f_2\chi_2+f_3\chi_3:\quad\supp f_k\subseteq l_{\zeta_k}\,(1\leq k\leq 3);\\
{\displaystyle \left.\|f\|^2=\sum\limits_k\int\limits_0^\infty|f_k(x\zeta_k)|^2dx<\infty\right\}}
\end{array}\label{eq1.14}
\end{equation}
with scalar product
$$\langle f,g\rangle=\sum\limits_k\int\limits_0^\infty f_k(x\zeta_k)\overline{g_k(x\zeta_k)}dx.$$
Match every component $f_k$ ($1\leq k\leq3$) of the element $f\in L^2(l)$ with its Fourier transform,
\begin{equation}
\widetilde{f}_k(\lambda)=(\mathcal{F}_{\zeta_k}f_k)(\lambda)\stackrel{\rm def}{=}\int\limits_0^\infty e^{-i\lambda\zeta_kx}f_k(x\zeta_k)dx\quad(1\leq k\leq3)\label{eq1.15}
\end{equation}

Description of $\{\widetilde{f}_k(\lambda)\}$ \eqref{eq1.15} is as follows.

(A) The function $f_1(\lambda)$ \eqref{eq1.15} is holomorphically extendable into the lower half-plane,
\begin{equation}
\mathbb{C}_-(\zeta_1)(=\mathbb{C}_-)\stackrel{\rm def}{=}\{\lambda=\mu+i\nu\in\mathbb{C}:\nu<0\},\label{eq1.16}
\end{equation}
and due to Paley -- Wiener theorem \cite{22, 23}, the operator $\mathcal{F}_{\zeta_1}$ \eqref{eq1.15} is a unitary isomorphism between $L^2(l_{\zeta_1})$ and the space $H_-^2(\zeta_1)(=H_-^2)$, besides, $\|f_1\|_{L^2(l_{\zeta_1})}=2\pi\|\widetilde{f}_1\|^2_{L^2(L_{\zeta_1})}$. The inverse operator $\mathcal{F}_{\zeta_1}^{-1}$ is given by
$$(\mathcal{F}_{\zeta_1}^{-1}g)(x)=\frac1{2\pi}\int\limits_{L_{\zeta_1}}e^{i\lambda\zeta_1x}g(\lambda\zeta_1)d\lambda\quad(g\in H(\zeta_1))$$
and $\mathcal{F}_{\zeta_1}^{-1}=(\mathcal{F}_{\zeta_1})^*$.

(B) The number $\lambda\zeta_2$ is real only if $\lambda\in L_{\zeta_3}$ \eqref{eq1.9}, i. e., $\lambda=\zeta_3\eta$ ($\eta\in\mathbb{R}$) and $\lambda\zeta_2=\eta\in\mathbb{R}$, therefore
$$\widetilde{f}_2(\lambda)=\int\limits_0^\infty e^{-i\eta x}f_2(x\zeta_2)dx\quad(\lambda=\zeta_3\eta\in L_{\zeta_3}),$$
and thus the function $\widetilde{f}_2(\lambda)$ ($\lambda=\zeta_3\eta\in L_{\zeta_3}$) on the straight line $L_{\zeta_3}$ is an ordinary Fourier transform. The function $\widetilde{f}_2(\lambda)$ has analytic extension into the half-plane
\begin{equation}
\mathbb{C}_-(\zeta_3)\stackrel{\rm def}{=}\{\lambda=\zeta_3\eta:\eta=\mu+i\nu,\nu<0\}\label{eq1.17}
\end{equation}
over the straight line $L_{\zeta_3}$, besides, $\lambda\in\mathbb{C}_-(\zeta_3)\Longleftrightarrow\eta\in\mathbb{C}_-(\zeta_1)$ \eqref{eq1.16} where $\lambda=\zeta_3\eta$. The straight line
$$L_{\zeta_3}(a)\stackrel{\rm def}{=}\{\lambda=\zeta_3(\eta-ia):\eta,a\in\mathbb{R}\}$$
is parallel to $L_{\zeta_3}$ and for $a>0$ it lies in $\mathbb{C}_-(\zeta_3)$ \eqref{eq1.17}, besides,
$$\widetilde{f}_2(\lambda)=\int\limits_0^\infty e^{-i\eta x}e^{-ax}f_2(x\zeta_2)dx\quad(\forall\lambda\in L_{\zeta_3}(a)).$$
Therefore integrals of the squares of moduli of $\widetilde{f}_2(\lambda)$ along each straight line $L_{\zeta_3}(a)$ ($a>0$) are uniformly bounded. Such functions form a Hardy class $H_-^2(\zeta_3)\subset L^2(L_{\zeta_2})$ corresponding to the half-plane $\mathbb{C}_-(\zeta_3)$. So, the operator $\mathcal{F}_{\zeta_2}$ \eqref{eq1.15} is an isomorphism between spaces $H^2(l_{\zeta_2})$ and $H_-^2(\zeta_3)$, besides,
$$(\mathcal{F}_{\zeta_2}^{-1}g)(x)=\frac1{2\pi}\int\limits_{L_{\zeta_3}}e^{i\lambda\zeta_2x}g(\lambda\zeta_3)d\lambda\quad(g\in H_-^2(\zeta_3)),$$
and Parseval's identity becomes
$$\langle\mathcal{F}_{\zeta_2}(f_2)(\lambda),\mathcal{F}_{\zeta_2}(\widetilde{f}_2(\lambda))\rangle_{L^2(L_{\zeta_3})}=2\pi\langle f_2(x),\widetilde{f}_2(x)\rangle_{L^2(l_{\zeta_2})}$$
($\forall f_2$, $\widetilde{f}_2\in L^2(l_{\zeta_2}))$).

(C) The number $\lambda\zeta_3$ is real if $\lambda\in L_{\zeta_2}$ \eqref{eq1.8}, i. e., $\lambda=\zeta_2\eta$ ($\eta\in\mathbb{R}$) and $\lambda\zeta_3=\eta\in\mathbb{R}$, thus
$$\widetilde{f}_3(\lambda)=\int_0^\infty e^{-i\eta x}f_3(x\zeta_3)dx\quad(\lambda=\zeta_2\eta\in L_{\zeta_2}).$$
Hence it follows that $\widetilde{f}_3(\lambda)$ is holomorphically extendable into the half-plane
\begin{equation}
\mathbb{C}_-(\zeta_2)\stackrel{\rm def}{=}\{\lambda=\zeta_2\eta:\eta=\mu+i\nu,\nu<0\}\label{eq1.18}
\end{equation}
over the straight line $L_{\zeta_2}$ and $\lambda\in\mathbb{C}_-(\zeta_2)\Longleftrightarrow\eta\in\mathbb{C}_-(\zeta_1)$ \eqref{eq1.16} ($\lambda=\zeta_2\eta$). The function $\widetilde{f}_3(\lambda)$ is of a Hardy space $H_-^2(\zeta_2)\subset L^2(L_{\zeta_2})$ corresponding to the half-plane $\mathbb{C}_-(\zeta_2)$ \eqref{eq1.18}. As in (B), the operator $\mathcal{F}_{\zeta_3}$ \eqref{eq1.15} establishes isomorphism between the spaces $L^2(L_{\zeta_3})$ and $H_-^2(\zeta_2)$ and the inverse to it operator is
$$(\mathcal{F}_{\zeta_3}^{-1}g)(x)=\frac1{2\pi}\int\limits_{L_{\zeta_2}}e^{i\lambda\zeta_3x}g(\lambda\zeta_2)d\lambda\quad(g\in H_-^2(\zeta_2)),$$
besides, Parseval's identity is

\begin{picture}(200,200)
\put(0,100){\vector(1,0){200}}
\put(110,100){\vector(1,0){10}}
\put(150,0){\vector(-1,2){100}}
\put(100,100){\vector(-1,2){10}}
\put(100,100){\vector(-1,-2){10}}
\put(150,200){\vector(-1,-2){100}}
\put(30,190){$L_{\zeta_2}$}
\put(65,190){$\mathbb{C}_-(\zeta_2)$}
\put(25,160){$\mathbb{C}_+(\zeta_2)$}
\put(80,110){$\zeta_2$}
\put(90,70){$\zeta_3$}
\put(115,89){$\zeta_1$}
\put(190,105){$L_{\zeta_1}$}
\put(140,105){$\mathbb{C}_+(\zeta_1)$}
\put(140,85){$\mathbb{C}_-(\zeta_1)$}
\put(20,30){$\mathbb{C}_-(\zeta_3)$}
\put(70,30){$\mathbb{C}_+(\zeta_3)$}
\put(60,0){$L_{\zeta_3}$}
\end{picture}

\hspace{20mm}Fig. 1.

$$\langle(\mathcal{F}_{\zeta_3}f_3)(\lambda),(\mathcal{F}_{\zeta_3}\widehat{f}_3)(\lambda)\rangle_{L^2(L_{\zeta_2})}=2\pi\langle f_3(x),\widehat{f}_3(x)\rangle_{L^2(l_{\zeta_3})}.$$
Define the Hilbert space
\begin{equation}
\begin{array}{ccc}
L^2(L)\stackrel{\rm def}{=}\{\widetilde{f}=\sum\limits_k\widetilde{f}_k\widetilde{\chi}_k:\supp\widetilde{f}_k\in L_{\zeta_k}(1\leq k\leq3);\,\\
{\displaystyle\left.\|f\|^2=\sum\limits_k\int\limits_{L_{\zeta_k}}|\widetilde{f}_k(x\zeta_k)|^2dx<\infty\right\}}
\end{array}\label{eq1.19}
\end{equation}
given on the bundle of straight lines
\begin{equation}
L\stackrel{\rm def}{=}\bigcup\limits_k L_{\zeta_k}\label{eq1.20}
\end{equation}
where $\{L_{\zeta_k}\}$ are given by \eqref{eq1.8}, $\widetilde{\chi}_k=\chi_{L_{\zeta_k}}(x)$ are characteristic functions of the sets $L_{\zeta_k}$ ($1\leq k\leq3$).

\begin{theorem}\label{t1.1}
Fourier transform $\{\mathcal{F}_{\zeta_k}\}$ \eqref{eq1.15} is a unitary isomorphism of the space $L^2(l)$ \eqref{eq1.14} onto the subspace $H_-^2(\zeta_1)\oplus H_-^2(\zeta_2)\oplus H_-^2(\zeta_3)$ in $L^2(L)$ \eqref{eq1.19}.
\end{theorem}

For the space $L^2(\widehat{l})$, where
\begin{equation}
\widehat{l}\stackrel{\rm}{=}\bigcup\limits_k\widehat{l}_{\zeta_k}\label{eq1.21}
\end{equation}
($\widehat{l}_{\zeta_k}$ are given by \eqref{eq1.9}), analogue of this theorem is true with natural substitution of Hardy spaces $H_-^2(\zeta_k)\rightarrow H_+^2(\zeta_k)$ ($1\leq k\leq3$). Theorem \ref{t1.1} is a generalization of the Paley -- Wiener theorem onto the space of functions given on the bundle of rays $l$ \eqref{eq1.13}.

\begin{remark}\label{r1.2}
Since $\bigcap\limits_k\mathbb{C}_-(\zeta_k)=\{0\}$, where $\mathbb{C}_-(\zeta_k)$ are given by \eqref{eq1.16} -- \eqref{eq1.18}, Fourier transforms $\{\widetilde{f}_k(\lambda)\}$ \eqref{eq1.15} simultaneously exist only for $\lambda=0$.
\end{remark}
\vspace{5mm}

{\bf 1.3} Define the space
\begin{equation}
\begin{array}{ccc}
{\displaystyle L^2(l,a)\stackrel{\rm def}{=}\left\{f=\sum\limits_kf_k\chi_k:\supp f_k\subseteq l_{\zeta_k}(1\leq k\leq3);\right.}\\
{\displaystyle\left.\sum\limits_k\int\limits_0^\infty e^{2ax}|f_k(x\zeta_k)|^2dx<\infty\right\}}
\end{array}\label{eq1.22}
\end{equation}
where $a\in\mathbb{R}$. For $b\geq a$, $L^2(l,b)\subseteq L^2(l,a)$, particularly, $L^2(l,a)\subseteq L^2(l)$, for all $a\geq0$. Consider the Fourier transform \eqref{eq1.15} of the components of a function $f$ from $L^2(l,a)$ considering $a\geq0$.

($A_1$) The equality
$$\widetilde{f}_1(\lambda)=\int\limits_0^\infty e^{-i\lambda\zeta_1x}f_1(x\zeta_1)dx=\int\limits_0^\infty e^{-i\zeta_1(\lambda-ia)x}e^{ax}f_1(x\zeta_1)dx$$
and $e^{ax}f_1(x\zeta_1)\in L^2(l_{\zeta_1})$ imply that $\widetilde{f}_1(\lambda)$ is holomorphically extendable into the half-plane
\begin{equation}
\mathbb{C}_-(\zeta_1,a)\stackrel{\rm def}{=}\{\lambda=\mu+i\nu\in\mathbb{C}:\nu<a\}\label{eq1.23}
\end{equation}
corresponding to the Hardy class $H_-^2(\zeta_1,a)$ \cite{23}.

($B_1$) For $\widetilde{f}_2(\lambda)$ \eqref{eq1.15},
$$\widetilde{f}_2(\lambda)=\int\limits_0^\infty e^{-i\lambda\zeta_2x}f_2(x\zeta_2)dx=\int\limits_0^\infty e^{-i\zeta_2(\lambda-i\zeta_3a)}e^{ax}f_2(x\zeta_2)dx,$$
therefore $e^{ax}f(x\zeta_2)\in L^2(l_{\zeta_2})$ yields (see (B)) that $\widetilde{f}_2(\lambda)$ is analytically extendable into the half-plane
\begin{equation}
\mathbb{C}_-(\zeta_2,a)\stackrel{\rm def}{=}\{\lambda=\zeta_3\eta:\eta=\mu+i\nu,\nu<a\},\label{eq1.24}
\end{equation}
and thus $\widetilde{f}_2(\lambda)$ is of the Hardy class $H_-^2(\zeta_3,a)$ corresponding to $\mathbb{C}_-(\zeta_3,a)$ \eqref{eq1.24}.

($C_1$) Analogously, $\widetilde{f}_3(\lambda)$ \eqref{eq1.15}, for $f\in L^2(l,a)$, has the holomorphic extension into the half-plane
\begin{equation}
\mathbb{C}_-(\zeta_2,a)\stackrel{\rm def}{=}\{\lambda=\zeta_2\eta:\eta=\mu+i\nu,\nu<a\}\label{eq1.25}
\end{equation}
and is of the Hardy class $H_-^2(\zeta_2,a)$ which corresponds to $\mathbb{C}_-(\zeta_2,a)$ \eqref{eq1.25}.
Intersection of the half-planes $\bigcap \mathbb{C}_-(\zeta_k,a)$ ($\mathbb{C}_-(\zeta_k,a)$ are given by \eqref{eq1.23} -- \eqref{eq1.25}) produces an equilateral triangle $T_a$ (see Fig. 2) with side $2\sqrt3a$,
\begin{equation}
T_a\stackrel{\rm def}{=}\bigcap\limits_k\overline{C_-(\zeta_k,a)}=\{\lambda=\mu+i\nu\in\mathbb{C}:\nu\leq a,\nu\geq\sqrt3\mu-2a,\nu\geq-\sqrt\mu-2a\}.\label{eq1.26}
\end{equation}
As a result, we obtain the following statement.

\begin{picture}(200,200)
\put(0,100){\vector(1,0){200}}
\put(0,130){\line(1,0){200}}
\put(150,0){\vector(-1,2){100}}
\put(120,0){\line(-1,2){100}}
\put(180,200){\line(-1,-2){100}}
\put(150,200){\vector(-1,-2){100}}
\put(30,190){$L_{\zeta_2}$}
\qbezier[90](100,0)(100,100)(100,200)
\put(125,190){$\mathbb{C}_-(\zeta_3,a)$}
\put(40,160){$\mathbb{C}_-(\zeta_2,a)$}
\put(190,105){$L_{\zeta_1}$}
\put(142,120){$\mathbb{C}_-(\zeta_1,a)$}
\put(102,130){$a$}
\put(100,40){$-2a$}
\put(110,105){$T_a$}
\put(60,0){$L_{\zeta_3}$}
\end{picture}

\hspace{20mm} Fig. 2.

\begin{theorem}\label{t1.2}
The Fourier transform $\{\mathcal{F}_{\zeta_k}\}$ \eqref{eq1.15} sets the isomorphism of $L^2(l,a)$ onto $H_-^2(\zeta_1,a)\oplus H_-^2(\zeta_2,a)\oplus H_-^2(\zeta_3,a)$ in $L^2(L)$ \eqref{eq1.19}. At every point $\lambda$ of the triangle $T_a$ \eqref{eq1.26}, all Fourier transforms $\{\widetilde{f}_k(\lambda)\}$ \eqref{eq1.15} simultaneously exist and $\widetilde{f}_k(\lambda)$ ($1\leq k\leq3$) are holomorphic at every inner point $\lambda\in T_a$.
\end{theorem}
\vspace{5mm}

{\bf 1.4} Let $\widetilde{f}_1(\lambda)$ be Fourier transform \eqref{eq1.15} of the first component of $f\in L^2(l,a)$ \eqref{eq1.22} ($a\geq0$), then the functions $\widetilde{f}_1(-\lambda)$ form the Hardy space $H_+^2(\zeta_1,-a)$ corresponding to the half-plane
\begin{equation}
\mathbb{C}_+\stackrel{\rm def}{=}\{\lambda=\mu+i\nu\in\mathbb{C}:\nu>-a\}.\label{eq1.27}
\end{equation}
Analogously, the functions $\widetilde{f}_2(-\lambda)$ ($\widetilde{f}_2(\lambda)$ is Fourier transform \eqref{eq1.15} of the second component of $f\in L^2(l,a)$) set the Hardy class $H_+^2(\zeta_3,-a)$ corresponding to the half-plane
\begin{equation}
\mathbb{C}_+(\zeta_3,-a)\stackrel{\rm def}{=}\{\lambda=\zeta_3\eta:\eta=\mu+i\nu;\nu>-a\}.\label{eq1.28}
\end{equation}
Finally, $\widetilde{f}_3(-\lambda)$ ($\widetilde{f}_3(\lambda$ is Fourier transform \eqref{eq1.15} of the third component of $f\in L^2(l,a)$) generate the Hardy space $H_+^2(\zeta_2,-a)$ corresponding to the half-plane
\begin{equation}
\mathbb{C}_+(\zeta_2,-a)\stackrel{\rm def}{=}\{\lambda=\zeta_2\eta:\eta=\mu+i\nu,\nu>-a\}.\label{eq1.29}
\end{equation}
\begin{picture}(200,200)
\put(0,100){\vector(1,0){200}}
\put(0,70){\line(1,0){200}}
\put(150,0){\vector(-1,2){100}}
\put(180,0){\line(-1,2){100}}
\put(120,200){\line(-1,-2){100}}
\put(150,200){\vector(-1,-2){100}}
\put(30,190){$L_{\zeta_2}$}
\qbezier[90](100,0)(100,100)(100,200)
\put(52,190){$\mathbb{C}_+(\zeta_2,-a)$}
\put(8,43){$\mathbb{C}_+(\zeta_3,-a)$}
\put(190,105){$L_{\zeta_1}$}
\put(142,80){$\mathbb{C}_+(\zeta_1,-a)$}
\put(100,65){$-a$}
\put(75,80){$T_a^*$}
\put(60,0){$L_{\zeta_3}$}
\end{picture}

\hspace{20mm} Fig. 3

Intersection of the half-planes $\{\mathbb{C}_+(\zeta_k,-a)\}$ \eqref{eq1.27} -- \eqref{eq1.29} gives an equilateral triangle $T_a^*$ (see Fig. 3) with the side $2\sqrt3a$,
\begin{equation}
T_a^*\stackrel{\rm def}{=}\bigcap\limits_k\overline{\mathbb{C}_+(\zeta_k,-a)}=\{\lambda=\mu+i\nu\in\mathbb{C}:\nu\geq-a,\nu\leq\sqrt3\mu-2a,\nu\leq\sqrt3\mu+2a\}.\label{eq1.30}
\end{equation}
The triangle $T_a^*$ is complexly adjoint to $T_a$ \eqref{eq1.26}; $\lambda\in T_a\Longleftrightarrow\overline{\lambda}\in T_a^*$.

\begin{theorem}\label{t1.3}
For every function $f\in L^2(l,a)$ \eqref{eq1.23} ($a\geq0$) for all $\lambda\in T_a\cap T_a^*$ ($T_a$ and $T_a^*$ are given by \eqref{eq1.26} and \eqref{eq1.30} correspondingly), there exist Fourier transforms $\{\widetilde{f}_k(\lambda)\}$ \eqref{eq1.15} and $\{\widetilde{f}_k(-\lambda)\}$ which are analytic functions
at all inner points from $T_a\cap T_a^*$.
\end{theorem}

Intersection $T_a\cap T_a^*$ is a regular hexagon with the side $2a/\sqrt3$ (see Figs 2, 3).

\begin{remark}\label{r1.3}
Instead of $L^2(l,a)$ \eqref{eq1.22}, one can consider the space
\begin{equation}
L^2(l,a_k)\stackrel{\rm def}{=}\left\{f=\sum f_k\chi_k(\supp f_k\subseteq l_{\zeta_k}):\sum\limits_k\int\limits_0^\infty e^{2a_kx}|f_k(x_k\zeta_k)|^2dx<\infty\right\}\label{eq1.31}
\end{equation}
where $a_k\in\mathbb{R}$ ($1\leq k\leq3$). For $a_k\geq0$ ($\forall k$) analogues of Theorems \ref{t1.2}, \ref{t1.3} hold, but in this case the triangle $T_a$ \eqref{eq1.26} no longer is equilateral.
\end{remark}

\section{Jost solutions}\label{s2}

{\bf 2.1} In the Hilbert space
\begin{equation}
\mathcal{H}\stackrel{\rm def}{=}L^2(\mathbb{R}_-)\oplus L^2(\mathbb{R}_+)=\{y=(v,u):v\in L^2(\mathbb{R}_-),u\in L^2(\mathbb{R}_+)\},\label{eq2.1}
\end{equation}
consider the linear operator
\begin{equation}
\mathcal{L}_qy=(-iDv,iD^3u+qu)\label{eq2.2}
\end{equation}
where $y=(v,u)\in\mathcal{H}$; $D=d/dx$; $q=q(x)$ is a real function satisfying the condition
\begin{equation}
\int\limits_{\mathbb{R}_+}|q(x)|^2e^{2ax}dx<\infty\label{eq2.3}
\end{equation}
for some $a\geq0$. Domain $\mathfrak{D}(\mathcal{L}_q)$ of the operator $\mathcal{L}_q$ is given by the boundary conditions,
\begin{equation}
\mathfrak{D}(\mathcal{L}_q)\stackrel{\rm def}{=}\{y=(v,u)\in\mathcal{H}:v\in W_2^1(\mathbb{R}_-);u\in W_2^3(\mathbb{R}_+);u(0)=0,u'(0)=\theta v(0)\}\label{eq2.4}
\end{equation}
where $\theta\in\mathbb{T}$. Operator $\mathcal{L}_q$ is self-adjoint.

\begin{remark}\label{r2.1}
Relation \eqref{eq2.3} implies that $q(x)e^{bx}\in L^2(\mathbb{R}_+)$ for all such $b$ that $0\leq b\leq a$. Moreover, for all $c$ ($0\leq c<a$), $e^{cx}q(x)\in L^1(\mathbb{R}_+)$, because
$$\int\limits_{\mathbb{R}_+}e^{cx}|q(x)|dx=\int\limits_{\mathbb{R}_+}e^{ax}\cdot q(x)e^{(c-a)x}dx<\infty,$$
due to the Cauchy -- Bunyakovsky -- Schwarz inequality, since $e^{ax}|q(x)|$ and $e^{(c-a)x}$ belong to $L^2(\mathbb{R}_+)$.
\end{remark}

\begin{remark}\label{r2.2}
Symmetric operator $iD^3+q$ in $L^2(\mathbb{R}_+)$, domain of which consists of functions $u\in W_2^3(\mathbb{R}_+)$ such that $u(0)=u'(0)=0$, does not have self-adjoint extensions in $L^2(\mathbb{R}_+)$ \cite{16, 17, 3, 18, 19, 21}. Such extensions are constructed with exit from the space $L^2(\mathbb{R}_+)$, and the operator $\mathcal{L}_q$ \eqref{eq2.2}, \eqref{eq2.4} is one of such extensions.

General construction of self-adjoint extensions of operators $\mathcal{L}$ in $H$ is as follows. In the sum of Hilbert spaces $E\oplus H$, we consider operator $\mathcal{A}f=Ag+\mathcal{L}h$ ($f=(g,h)\in E\oplus H$). We select domain of the operator $A$ such that $\langle\mathcal{A}f,f'\rangle-\langle f,\mathcal{A}f'\rangle=0$. To do this, we have to 'match' \cite{3, 16, 17, 18, 19, 21} the boundary form $\langle\mathcal{L}h,h'\rangle-\langle h,\mathcal{L}h'\rangle$ of the operator $\mathcal{L}$ in $H$ with the boundary form $\langle Ag,g'\rangle-\langle g,Ag'\rangle$ of the operator $A$ in $E$. This 'matching' is characterized by the choice of the {\bf boundary operator} $B$ which maps '{\bf boundary defect elements}' of the operator $\mathcal{L}$ onto 'boundary elements' of the operator $A$ such that difference of the boundary forms $\mathcal{L}$ and $A$ vanishes. Hence, the self-adjoint extensions of the operator $\mathcal{L}$ are parameterized by the pair: by the operator $A$ in $H$ and boundary operator of matching $B$ acting in the space of boundary elements.
\end{remark}

Consider the equation system in $\mathcal{H}$,
\begin{equation}
\left\{
\begin{array}{lll}
-iDv=\lambda^3v\quad(x\in\mathbb{R}_-);\\
iD^3u+qu=\lambda^3u\quad(x\in\mathbb{R}_+).
\end{array}\right.\label{eq2.5}
\end{equation}
Solution to the first equation in \eqref{eq2.5} is given by $ce^{i\lambda^3x}$ ($c\in\mathbb{C}$, $x\in\mathbb{R}_-$), and the second equation, for $x\rightarrow\infty$ (due to \eqref{eq2.3}), transforms into equation \eqref{eq1.1}. Therefore asymptotic behavior of the solution $u(x)$, when $x\rightarrow\infty$, is described by the linear combination of the exponents $\{e^{i\lambda\zeta_0x}\}$.

By $e_k(\lambda,x)$, we denote the {\bf Jost solutions} \cite{1, 2, 3, 4, 5, 6, 19} to the equation
\begin{equation}
iD^3u(x)+q(x)u(x)=\lambda^3u(x)\quad(x\in\mathbb{R}_+,\lambda\in\mathbb{C})\label{eq2.6}
\end{equation}
satisfying the boundary condition
\begin{equation}
\lim\limits_{x\rightarrow\infty}(u(x)-e^{i\lambda\zeta_px})=0\quad(1\leq p\leq3).\label{eq2.7}
\end{equation}
Show that the integral equation
\begin{equation}
e_p(\lambda,x)=e^{i\lambda\zeta_px}-i\int\limits_x^\infty\frac{s_z(i\lambda(x-t))}{(i\lambda)^2}q(t)e_p(\lambda,t)dt\quad(1\leq p\leq3)\label{eq2.8}
\end{equation}
is equivalent to boundary value problem \eqref{eq2.6}, \eqref{eq2.7}. First, we establish solvability of equation \eqref{eq2.8}. To do this, in the space $L^2(\mathbb{R}_+)$ set a family of Volterra operators (depending on $\lambda\in\mathbb{C}$)
\begin{equation}
(K_\lambda f)(x)\stackrel{\rm def}{=}\int\limits_x^\infty K_1(\lambda,x,t)q(t)f(t)dt\quad(f\in L^2(\mathbb{R}_+))\label{eq2.9}
\end{equation}
where
\begin{equation}
K_1(\lambda,x,t)\stackrel{\rm def}{=}\frac{s_2(i\lambda(x-t))}{(i\lambda)^2}.\label{eq2.10}
\end{equation}
Equation \eqref{eq2.8}, in terms of $K_\lambda$ \eqref{eq2.9}, becomes
$$(I+K_\lambda)e_p(\lambda,x)=e^{i\lambda\zeta_px},$$
and thus
\begin{equation}
e_p(\lambda,x)=\sum\limits_0^\infty(-i)^nK_\lambda^ne^{i\lambda\zeta_px}.\label{eq2.11}
\end{equation}
The operators $K_\lambda^n$ also are Volterra,
$$(K_\lambda^nf)(x)\stackrel{\rm def}{=}\int\limits_x^\infty K_n(\lambda,x,t)q(t)f(t)dt\quad(f\in L^2(\mathbb{R}_+)),$$
and for the kernels $K_n(\lambda,x,t)$ recurrent relations take place,
\begin{equation}
K_{n+1}(\lambda,x,t)=\int\limits_x^tK_n(\lambda,x,s)q(s)K_1(\lambda,s,t)ds\quad(n\in\mathbb{N}).\label{eq2.12}
\end{equation}
Since
\begin{equation}
\begin{array}{ccc}
i\lambda\zeta_1=-\beta+i\alpha;\, i\lambda\zeta_2=\frac12[(\beta-\alpha\sqrt3)-i(\alpha+\beta\sqrt3)];\\ i\lambda\zeta_3=\frac12[(\beta+\alpha\sqrt3)-i(\alpha-\beta\sqrt3)]
\end{array}\label{eq2.13}
\end{equation}
($\lambda=\alpha+i\beta$; $\alpha$, $\beta\in\mathbb{R}$),
$$s_p(i\lambda)=\frac13\left\{e^{-\beta+\alpha i}+\frac1{\zeta_2^p}e^{\frac12((\beta-\alpha\sqrt3)-i(\alpha-\beta\sqrt3))}+\frac1{\zeta_3^p}e^{\frac12((\beta+\alpha\sqrt3)-i(\alpha+\beta\sqrt3))}\right\},$$
hence
$$|s_p(i\lambda)|\leq\frac13\left(e^{-\beta}+2e^{\frac\beta2}\ch\frac{\alpha\sqrt3}2\right),$$
therefore
\begin{equation}
|s_p(i\lambda)|\leq d(\lambda)\quad\left(d(\lambda)\stackrel{\rm def}{=}e^{|\beta|}\ch\frac{\alpha\sqrt3}2;\,\lambda=\alpha+i\beta\right).\label{eq2.14}
\end{equation}
So, for the kernel $K_1(\lambda,x,t)$ \eqref{eq2.10}, the following estimation holds,
\begin{equation}
|K_1(\lambda,x,t)|\leq\frac{d(\lambda(t-x))}{|\lambda|^2}\quad(\lambda\not=0,t\geq x).\label{eq2.15}
\end{equation}

\begin{lemma}\label{l2.1}
Kernels $K_n(\lambda,x,t)$ \eqref{eq2.12} have the properties
\begin{equation}
K_n(\lambda\zeta_2,x,t)=K_n(\lambda,x,t);\quad\overline{K_n(\lambda,x,t)}=K_n(\overline{\lambda},t,x)\label{eq2.16}
\end{equation}
and satisfy the inequalities
\begin{equation}
|K_n(\lambda,x,t)|\leq\left\{
\begin{array}{lll}
{\displaystyle\frac{d(\lambda(t-x))}{|\lambda|^{2n}}\cdot\frac{\sigma^{n-1}(t)}{(n-1)!}\quad(\lambda\not=0,t\geq x);}\\
{\displaystyle\left(\frac{(t-x)^2}2\right)^n\frac{\sigma^{n-1}(t)}{n^{2n}(n-1)!}\quad(\lambda=0)}
\end{array}\right.\label{eq2.17}
\end{equation}
($n\in\mathbb{N}$) where $d(\lambda)$ is given by \eqref{eq2.14} and
\begin{equation}
\sigma(t)\stackrel{\rm def}{=}\int\limits_0^t|q(s)|ds.\label{eq2.18}
\end{equation}
\end{lemma}

P r o o f. The kernel $K_1(\lambda,x,t)$ \eqref{eq2.10}, due to (ii), (iii) \eqref{eq1.4}, satisfy the equalities \eqref{eq2.16}, hence, in view of \eqref{eq2.12}, it follows that $K_n(\lambda,x,t)$ also has properties \eqref{eq2.16} ($q(x)$ is real).

Let $\lambda\not=0$, then inequality \eqref{eq2.17}, for $n=1$, coincides with \eqref{eq2.15}. Using induction by $n$ and recurrent relation \eqref{eq2.12}, and estimates \eqref{eq2.15}, \eqref{eq2.17} (for $n$), obtain
$$|K_{n+1}(\lambda,x,t)|\leq\int\limits_x^t\frac{d(\lambda(s-1))}{|\lambda|^{2n}}\frac{\sigma^{n-1}(s)}{(n-1)!}|q(s)|\frac{d(\lambda(t-s))}{|\lambda|^2}ds.$$
Since
$$d(\lambda(s-x))d(\lambda(t-s))=e^{|\beta|(t-x)}\ch\frac{\alpha\sqrt3}2(s-x)\ch\frac{\alpha\sqrt3}2(t-s)$$
$$=\frac12e^{|\beta|(t-x)}\left\{\ch\frac{\alpha\sqrt3}2(t-x)+\ch\frac{\alpha\sqrt3}2(t+x-2s)\right\}\leq d(\lambda(t-x)),$$
then
$$|K_{n+1}(\lambda,x,t)|\leq\frac{d(\lambda(t-x))}{|\lambda|^{2(n+1)}}\cdot\frac{\sigma^n(t)}{n!},$$
and this gives \eqref{eq2.17} for $n+1$.

Consider the case of $\lambda=0$. Formula \eqref{eq2.10} implies that $K_1(0,x,t)=(t-x)^2/2$, this coincides with \eqref{eq2.17} ($n=1$, $\lambda=0$). Again using induction and \eqref{eq2.18}, we have
$$|K_{n+1}(0,x,t)|\leq\int\limits_x^t\left(\frac{(s-x)^2}2\right)^n\frac{\sigma^{n-1}(s)}{n^{2n}(n-1)!}\frac{(t-s)^2}2|q(s)|ds.$$
The function $f(s)=(t-s)^2(s-x)^{2n}$ is positive for $s\in(x,t)$ and $f(x)=f(t)=0$. It reaches its maximum when $s\in(x,t)$ at the point $s_0=(x+nt)/(n+1)$ and $f(s_0)=(t-x)^{2(n+1)}n^{2n}/(n+1)^{2(n+1)}$, therefore
$$|K_{n+1}(0,x,t)|\leq\left(\frac{(t-x)^2}2\right)^{n+1}\frac1{(n+1)^{2(n+1)}}\cdot\frac{\sigma^n(t)}{n!}.\blacksquare$$

Relation \eqref{eq2.11} implies that
\begin{equation}
e_p(\lambda,x)=e^{i\lambda\zeta_px}+\int\limits_x^\infty N(\lambda,x,t)q(t)e^{i\lambda\zeta_pt}dt\quad(1\leq p\leq3)\label{eq2.19}
\end{equation}
where
\begin{equation}
N(\lambda,x,t)=\sum\limits_{n=1}^\infty(-i)^nK_n(\lambda,x,t).\label{eq2.20}
\end{equation}
Series \eqref{eq2.20} is majorized by the converging series, due to \eqref{eq2.17}, and
\begin{equation}
|N(\lambda,x,t)|\leq\left\{
\begin{array}{lll}
{\displaystyle\frac{d(\lambda(t-x))}{|\lambda|^2}\cdot\exp\left\{\frac{\sigma(t)}{|\lambda|^2}\right\}\quad(\lambda\not=0);}\\
{\displaystyle\frac{(t-x)^2}2\sum\limits_0^\infty\left(\frac{(t-x)^2\sigma(t)}2\right)^n\frac1{(n+1)^{2(n+1)}n!}\quad(\lambda=0).}
\end{array}\right.\label{eq2.21}
\end{equation}
Hence it follows that integral \eqref{eq2.19} converges uniformly by $x$ when $|\lambda|<a$ (Remark \ref{r2.1}). Thus the solution $e_p(\lambda,x)$ to the integral equation \eqref{eq2.8} satisfies the differential equation \eqref{eq2.6}.

Turn to the boundary condition \eqref{eq2.7}. For $\lambda\not=0$, relations \eqref{eq2.19}, \eqref{eq2.21} imply that
$$|e_p(\lambda,x)|\leq e^{|\lambda|x}+\int\limits_x^\infty|q(t)|e^{|\lambda|t}\frac{d(\lambda(t-x))}{|\lambda|^2}\exp\left\{\frac{\sigma(t)}{|\lambda|^2}\right\}dt\leq e^{|\lambda|x}$$
$$+\frac1{|\lambda|^2}\int\limits_x^\infty|q(t)|e^{|\lambda|t}e^{2|\lambda|(t-x)}\exp\left\{\frac{\sigma(t)}{|\lambda|^2}\right\}dt\leq e^{|\lambda|x}\left(1+\frac{\exp\left(\frac\sigma{|\lambda|^2}\right)}{|\lambda|^2}\int\limits_x^\infty e^{3|\lambda|(t-x)}|q(t)|dt\right),$$
here $\sigma\stackrel{\rm def}{=}\sigma(\infty)$. Using $e^{3|\lambda|(t-x)}|q(t)|\leq e^{3|\lambda|t}|q(t)|\in L^1(\mathbb{R}_+)$ for $3|\lambda|<a$ (Remark \ref{r2.1}), obtain
$$|e_p(\lambda,x)|<e^{|\lambda|x}\left\{1+\frac{q_1}{|\lambda|^2}\exp\left(\frac\sigma{|\lambda|^2}\right)\right\}\quad\left(\forall x\in\mathbb{R}_+,|\lambda|<\frac a3\right),$$
besides, $q_1\stackrel{\rm def}{=}\|e^{3|\lambda|t}q(t)\|_{L^1(\mathbb{R}_+)}$. Using this inequality, we have
$$\left|\int\limits_x^\infty\frac{s_2(i\lambda(x-t))}{(i\lambda)^2}q(t)e_p(\lambda,t)dt\right|\leq\left(1+\frac{q_1}{|\lambda|^2}\exp\left(\frac\sigma{|\lambda|^2}\right)\right)
\int\limits_x^\infty\frac{d(\lambda(t-x))}{|\lambda|^2}|q(t)|e^{|\lambda|t}dt$$
$$\leq\frac1{|\lambda|^2}\left(1+\frac{q_1}{|\lambda|^2}\exp\left(\frac\sigma{|\lambda|^2}\right)\right)\int\limits_x^\infty e^{3|\lambda|t}|q(t)|dt.$$
The last integral tends to zero as $x\rightarrow\infty$ since $e^{3|\lambda|t}q(t)\in L^1(\mathbb{R}_+)$ when $3|\lambda|<a$ (Remark \ref{r2.1}), this ensures, due to \eqref{eq2.8}, the fulfillment of the boundary condition \eqref{eq2.7}.

Denote by $\mathbb{D}_{a/3}$ a circle
\begin{equation}
\mathbb{D}_{a/3}=\{\lambda\in\mathbb{C}:|\lambda|<a/3\}\quad(a>0)\label{eq2.22}
\end{equation}
lying inside the triangle $T_a$ \eqref{eq1.26}.

\begin{theorem}\label{t2.1}
If a real function $q(x)$ satisfies condition \eqref{eq2.3}, then for all $\lambda\in\mathbb{D}_{a/3}$ \eqref{eq2.22} the Jost solutions $\{e_p(\lambda,x)\}$ of the boundary value problem \eqref{eq2.6}, \eqref{eq2.7} exist and are given by \eqref{eq2.19}. The functions $\{e_p(\lambda,x)\}$ are holomorphic by $\lambda$ for $\lambda\in\mathbb{D}_{a/3}$, and
\begin{equation}
e_p(\lambda\zeta_2,x)=e_{p'}(\lambda,x)\quad(p'=(p+1)\mod3,1\leq p\leq3).\label{eq2.23}
\end{equation}
\end{theorem}

P r o o f. Prove that the boundary condition \eqref{eq2.7} also holds for $\lambda=0$. To do this, write equation \eqref{eq2.8} as
$$e_p(0,x)-1=-i\int_x^\infty\frac{(t-x)^2}2q(t)[e_p(0,t)-1]dt-i\int\limits_x^\infty\frac{(t-x)^2}2q(t)dt.$$
Solving this equation with respect to $e_p'(0,x)-1$ by the method of successive approximations and using
$$\left|\int\limits_x^\infty\frac{(t-x)^2}2q(t)dt\right|\leq\int\limits_x^\infty\frac{t^2}2|q(t)dt,$$
we obtain
$$|e_p(t,x)-1|\leq\int\limits_x^\infty\frac{t^2}2|q(t)|dt\cdot\exp\left\{\int\limits_x^\infty\frac{t^2}2|q(t)|dt\right\}.$$
Using the inequality ${\displaystyle\frac{x^2}2<e^x}$ ($\forall x>0$), we find that
$$\int\limits_x^\infty\frac{t^2}2|q(t)|dt=\frac1{b^2}\int\limits_x^\infty\frac{(tb)^2}2|q(t)|dt\leq\frac1{b^2}\int\limits_x^\infty e^{tb}|q(t)|dt\quad(0<b<a),$$
and due to Remark \ref{r2.1} this integral tends to zero when $x\rightarrow\infty$. $\blacksquare$

The {\bf Jost solution} of system \eqref{eq2.5} is given by
\begin{equation}
E_p(\lambda,x)=(a_p(\lambda)e^{i\lambda^3x},e_p(\lambda,x))\quad(1\leq p\leq3),\label{eq2.24}
\end{equation}
here $\lambda\in\mathbb{D}_{a/3}$ \eqref{eq2.22}; $\{e_p(\lambda,x)\}$ are given by the formulas \eqref{eq2.19}; $\{a_p(\lambda)\}$ are some functions of $\lambda$. If $a_p(\lambda\zeta_2)=a_{p'}(\lambda)$ ($p'=(p+1)\mod3$), then, taking \eqref{eq2.23} into account, we have
\begin{equation}
E_p(\lambda\zeta_2,x)=E_{p'}(\lambda,x)\quad(p'=(p+1)\mod3,1\leq p\leq3.\label{eq2.25}
\end{equation}
\vspace{5mm}

{\bf 2.2} Study properties of Jost solutions.

\begin{lemma}\label{l2.2}
For all $\lambda\in\mathbb{D}_{a/3}\backslash\{0\}$, functions $\{e_p(\lambda,x)\}$ are linearly independent.
\end{lemma}

P r o o f. Assuming the contrary, suppose that, for some $\lambda\in\mathbb{D}_{a/3}\backslash\{0\}$, there are numbers $\{\mu_p\}_1^3$ from $\mathbb{C}$ (${\displaystyle\sum|\mu_p|\not=0}$) such that ${\displaystyle\sum\mu_pe_p(\lambda,x)=0}$ ($\forall x\in\mathbb{R}_+$), then
$$\sum\mu_pe_p(\lambda,x)=0,\quad\sum\mu_pe'_p(\lambda,x)=0,\quad\sum\mu_pe''_p(\lambda,x)=0.$$
Determinant of this system
\begin{equation}
\Delta(\lambda,x)\stackrel{\rm def}{=}\det\left[
\begin{array}{ccc}
e_1(\lambda,x)&e_2(\lambda,x)&e_3(\lambda,x)\\
e'_1(\lambda,x)&e'_2(\lambda,x)&e'_3(\lambda,x)\\
e''_1(\lambda,x)&e''_2(\lambda,x)&e''_3(\lambda,x)
\end{array}\right]\label{eq2.26}
\end{equation}
does not depend on $x$ since $\Delta'(\lambda,x)=0$, because $\{e_p(\lambda,x)\}$ are solutions to equation \eqref{eq2.6}, therefore, using \eqref{eq2.7}, we have
$$\Delta(\lambda,x)=\Delta_0(\lambda,x)+o\left(\frac1x\right)$$
where
$$\Delta(\lambda,\infty)=\det\left[
\begin{array}{lll}
e^{i\lambda\zeta_1x}&e^{i\lambda\zeta_2x}&e^{i\lambda\zeta_3x}\\
i\lambda\zeta_1e^{i\lambda\zeta_1x}&i\lambda\zeta_2e^{i\lambda\zeta_2x}&i\lambda\zeta_3e^{i\lambda\zeta_3x}\\
(i\lambda\zeta_1)^2e^{i\lambda\zeta_1x}&(i\lambda\zeta_2)^2e^{i\lambda\zeta_2x}&(i\lambda\zeta_3)^2e^{i\lambda\zeta_3x}
\end{array}\right]$$
$$=(i\lambda)^3\det\left[
\begin{array}{ccc}
1&1&1\\
1&\zeta_2&\zeta_3\\
1&\zeta_3&\zeta_2
\end{array}\right]=-3\sqrt3\lambda^3.$$
Thus $\Delta(\lambda,x)=-3\sqrt3\lambda^3$, and so for all $\lambda\in\mathbb{D}_{a/3}\backslash\{0\}$, determinant $\Delta(\lambda,x)$ of our system is not zero and its solution $\{\mu_p\}$ is trivial, $\mu_1=\mu_2=\mu_3=0$. $\blacksquare$

Turn to the analytic (relative to variable $\lambda$) properties of solutions $\{e_p(\lambda,x)\}$. Define the functions
\begin{equation}
\psi_p(\lambda,x)\stackrel{\rm def}{=}e_p(\lambda,x)e^{-i\lambda\zeta_px}\quad(1\leq p\leq3)\label{eq2.27}
\end{equation}
which, due to \eqref{eq2.8}, are solutions to the integral equations
\begin{equation}
\psi_p(\lambda,x)=1-i\int\limits_x^\infty e^{i\lambda\zeta_p(t-x)}\frac{s_2(i\lambda(x-t))}{(i\lambda)^2}q(t)\psi_p(\lambda,t)dt\quad(1\leq p\leq3).\label{eq2.28}
\end{equation}
Consider $\psi_1(\lambda,x)$, then kernel of equation \eqref{eq2.28}, for $p=1$, is
\begin{equation}
e^{i\lambda(t-x)}\frac{s_2(i\lambda(x-t))}{(i\lambda)^2}=\frac1{3(i\lambda)^2}\left\{1+\zeta_2e^{i\lambda(\zeta_1-\zeta_2)(t-x)}+\zeta_3e^{i\lambda(\zeta_1-\zeta_3)(t-x)}\right\}.
\label{eq2.29}
\end{equation}
Taking into account (see \eqref{eq1.2})
\begin{equation}
\frac{\zeta_2-\zeta_1}{\sqrt3}=-i\zeta_3;\quad\frac{\zeta_1-\zeta_3}{\sqrt3}=-i\zeta_2;\quad\frac{\zeta_3-\zeta_2}{\sqrt3}=-i\zeta_1\label{eq2.30}
\end{equation}
and $\lambda=\alpha+i\beta$ ($\alpha$, $\beta\in\mathbb{R}$), we obtain
\begin{equation}
\begin{array}{ccc}
{\displaystyle i\lambda(\zeta_1-\zeta_2)=-\sqrt3\lambda\zeta_3=\sqrt3\left[\frac12(\alpha-\sqrt3\beta)+\frac i2(\alpha\sqrt3+\beta)\right];}\\
{\displaystyle i\lambda(\zeta_1-\zeta_3)=\sqrt3\lambda\zeta_2=\sqrt3\left[\frac12(-\alpha-\sqrt3\beta)+\frac i2(\alpha\sqrt3-\beta)\right];}
\end{array}\label{eq2.31}
\end{equation}
hence it follows that
$$|e^{i\lambda(\zeta_1-\zeta_2)(t-x)}|\leq e^{\frac{\sqrt3}2(\alpha-\beta\sqrt3)(t-x)};\quad|e^{i\lambda(\zeta_1-\zeta_3)(t-x)}|\leq e^{\frac{\sqrt3}2(-\alpha-\beta\sqrt3)(t-x)}.$$
So, for all $\lambda=\alpha+i\beta\in\mathbb{C}$, such that $\alpha-\beta\sqrt3<0$ and $-\alpha-\beta\sqrt3<0$, each exponent in \eqref{eq2.29} is less than one in modulo (because $t-x\geq0$) and decreases exponentially as $\lambda\rightarrow\infty$ when $\lambda$ belongs to this set.

\begin{picture}(200,200)
\put(0,100){\vector(1,0){200}}
\put(100,200){\vector(0,-1){200}}
\put(200,67){\vector(-3,1){200}}
\put(0,67){\vector(3,1){200}}
\put(2,130){$-iL_{\zeta_3}$}
\put(100,150){$\psi_1(\lambda,x)$}
\put(70,150){$\Omega_2$}
\qbezier(153,120)(100,140)(49,120)
\qbezier(100,10)(170,20)(190,131)
\qbezier(45,120)(30,70)(100,25)
\qbezier(125,100)(115,110)(120,104)
\put(133,105){$\pi/6$}
\put(160,105){$S_2(i)$}
\put(20,50){$\Omega_3$}
\put(70,110){$S_4(i)$}
\put(70,70){$S_6(i)$}
\put(50,100){$S_5(i)$}
\put(0,80){$\psi_3(\lambda,x)$}
\put(115,70){$S_1(i)$}
\put(105,110){$S_3(i)$}
\put(180,135){$-iL_{\zeta_2}$}
\put(180,80){$\psi_2(\lambda,x)$}
\put(165,50){$\Omega_1$}
\put(60,0){$-iL_{\zeta_1}$}
\end{picture}

\hspace{20mm} Fig. 4.

Apply the method of successive approximations to equation \eqref{eq2.28} ($p=1$) and, using Lemma \ref{l2.1}, we obtain that $\psi_1(\lambda,x)$ is holomorphic in the domain $\alpha<\beta\sqrt3$, $-\alpha<\beta\sqrt3$. This set has a simple description. Denote by $S_p(i)$ the sectors obtained from $S_p$ \eqref{eq1.10} by `$-\pi/2$' rotation,
\begin{equation}
S_p(i)\stackrel{\rm def}{=}(-i)S_p\quad(1\leq p\leq6).\label{eq2.32}
\end{equation}
Upon this, straight lines $L_{\zeta_k}$ \eqref{eq1.8} transform into the straight lines $-iL_{\zeta_k}$ with directing vectors `$-i\zeta_k$'. By $\{\Omega_p\}$, we denote the sectors
\begin{equation}
\begin{array}{ccc}
\Omega_1\stackrel{\rm def}{=}S_1(i)\cup S_2(i)\cup(-i\widehat{l}_{\zeta_3});\quad\Omega_2\stackrel{\rm def}{=}S_3(i)\cup S_4(i)\cup (-i\widehat{l}_{\zeta_1});\\
\Omega_3\stackrel{\rm def}{=}S_5(i)\cup S_6(i)\cup(-i l_{\zeta_2}).
\end{array}\label{eq2.33}
\end{equation}

\begin{lemma}\label{l2.3}
For all $x\in\mathbb{R}_+$, the function $\psi_1(\lambda,x)$ \eqref{eq2.27} is analytic (relative to $\lambda$) in the sector $\Omega_2$ and the following asymptotic formulas hold,
\begin{equation}
\begin{array}{lll}
(a)&{\displaystyle\psi_1(\lambda,x)=1+\frac i{3\lambda^2}\int\limits_x^\infty q(t)dt+O\left(\frac1{|\lambda|^4}\right)}&(\lambda\rightarrow\infty,\lambda\in\overline{\Omega}_2);\\
(b)&{\displaystyle e'_1(\lambda,x)e^{-i\lambda\zeta_1x}=i\lambda\zeta_1-\frac1{3\lambda}\int\limits_x^\infty q(t)dt+O\left(\frac1{|\lambda|^3}\right)}&(\lambda\rightarrow\infty,\lambda\in\overline{\Omega}_2).
\end{array}\label{eq2.34}
\end{equation}
\end{lemma}

P r o o f. Using \eqref{eq2.29}, re-write equation \eqref{eq2.28} for $\psi_1(\lambda,x)$ as
$$\psi_1(\lambda,x)=1+\frac i{3\lambda^2}\int\limits_x^\infty\left\{1+\zeta_2e^{i\lambda(\zeta_1-\zeta_2)(t-x)}+\zeta_3e^{i\lambda(\zeta_1-\zeta_3)(t-x)}\right\}q(t)\psi_1(\lambda,t)dt.$$
For all $\lambda\in\Omega_2$, the exponents under the integral sign tend to zero when $\lambda\rightarrow\infty$ ($\lambda\in\Omega_2$), which proves relation (a) \eqref{eq2.34}. When $\lambda$ belongs to the border of $\partial\Omega_2$ of the sector $\Omega_2$, this equality follows from the Riemann -- Lebesgue lemma.

Since (see \eqref{eq2.8})
$$e'_1(\lambda,x)=i\lambda\zeta_1e^{i\lambda\zeta_1x}-i\int\limits_x^\infty\frac{s_1(i\lambda(x-t))}{i\lambda}q(t)e_1(\lambda,t)dt,$$
then
$$e'_1(\lambda,x)e^{-i\lambda\zeta_1x}=i\lambda\zeta_1-\frac1\lambda\int\limits_x^\infty e^{i\lambda(t-x)}s_1(i\lambda(x-t))\psi_1(\lambda,t)dt,$$
hence, due to analogous considerations, follows (b) \eqref{eq2.34}. $\blacksquare$

\begin{remark}\label{r2.3}
Equalities $\psi_1(\lambda\zeta_2,x)=\psi_2(\lambda,x)$, $\psi_1(\lambda\zeta_3,x)=\psi_3(\lambda,x)$, and \eqref{eq2.23} imply that $\psi_2(\lambda,x)$ and $\psi_3(\lambda,x)$ are holomorphic in $\Omega_1$ and $\Omega_3$ \eqref{eq2.33} correspondingly. Formulas analogous to \eqref{eq2.34}, for $e_2(\lambda,x)$ and $e_3(\lambda,x)$, are true in sectors $\overline{\Omega}_1$ and $\overline{\Omega}_3$.
\end{remark}

\begin{theorem}\label{t2.2}
The set of joint zeros of the functions $e_1(\lambda,0)$ and $e'_1(\lambda,0)$ lying inside the sector $\Omega_2$ \eqref{eq2.33} is finite and is given by
\begin{equation}
\Lambda_1\stackrel{\rm def}{=}\{\lambda_n,\mu_n\in\Omega_2:\lambda_n=-x_n\zeta_3\in\widehat{l}_{\zeta_3},\mu_n=x_n\zeta_2\in l_{\zeta_2}\},\label{eq2.35}
\end{equation}
where $l_{\zeta_2}$, $\widehat{l}_{\zeta_3}$ are rays \eqref{eq1.9} lying in the sector $\Omega_2$; $x_n$ are positive numbers enumerated in ascending order ($0\leq n\leq N<\infty$. Points $\{\lambda_n\}$, $\{\mu_n\}$ are simple zeros of the function $\theta_1(\lambda,0)$.
\end{theorem}

\begin{picture}(200,200)
\put(0,100){\vector(1,0){200}}
\put(100,0){\vector(0,1){200}}
\put(150,0){\vector(-1,2){100}}
\put(150,200){\line(-1,-2){100}}
\put(200,67){\line(-3,1){200}}
\put(0,67){\line(3,1){200}}
\put(30,190){$l_{\zeta_2}$}
\put(135,150){$\lambda_n$}
\put(135,130){$\Omega_2$}
\put(150,180){$\widehat{l}_{\zeta_3}$}
\put(121,150){$\times$}
\put(70,147){$\circ$}
\put(75,157){$\mu_n$}
\qbezier(153,120)(100,140)(49,120)
\end{picture}

\hspace{20mm} Fig. 5.

P r o o f. Analyticity of $e_1(\lambda,0)$ in the sector $\Omega_2$ (Lemma \ref{l2.3}) and $e_1(\lambda,0)=\psi_1(\lambda,0)\rightarrow1$ \eqref{eq2.34} as $\lambda\rightarrow\infty$ ($\lambda\in\Omega_2$) imply that the set of zeros of $e_1(\lambda,0)$ in $\Omega_2$ is bounded and no more than countable. Consider the Jost solution to system \eqref{eq2.5}
\begin{equation}
\varphi_1(\lambda,x)\stackrel{\rm def}{=}(v_1(\lambda,x),e_1(\lambda,x))\label{eq2.36}
\end{equation}
where $v_1(\lambda,x)=\overline{\theta}e'_1(\lambda,0)e^{i\lambda^3x}$ ($\theta\in\mathbb{T}$, $x\in\mathbb{R}_-$, $\lambda\in\Omega_2$). If $\lambda$ is a zero of the function $e_1(\lambda,0)$, then $\varphi_1(\lambda,x)$ \eqref{eq2.36} satisfies the boundary conditions \eqref{eq2.4} and
$$\mathcal{L}_q\varphi_1(\lambda,x)=\lambda^3\varphi_1(\lambda,x)\quad(\lambda\in\Omega_2).$$
The second component $e_1(\lambda,x)$ of the function $\varphi_1(\lambda,x)$ \eqref{eq2.36}, for $\lambda\in\Omega_2$ belongs to $L^2(\mathbb{R}_+)$, and the first component $v_1(\lambda,x)\in L^2(\mathbb{R}_-)$ if $\lambda\in\{\lambda\in\mathbb{C}:\pi/3<\arg\lambda<2\pi/3\}$ ($\subset\Omega_2$). At the points ${\displaystyle\lambda\in\Omega_2\backslash\left\{\lambda\in\mathbb{C}:\frac\pi3\leq\arg\lambda\leq\frac{2\pi}3\right\}}$, the function $v_1(\lambda,x)$ does not belong to $L^2(\mathbb{R}_-)$, excluding those $\lambda$ where $e'_1(\lambda,0)=0$, i. e., $v_1(\lambda,x)\equiv0$. Using the boundary conditions for $e_1(\lambda,0)$ and $e'_1(\lambda,0)=\theta v_1(\lambda,0)$, we obtain
$$(\lambda^3-\overline{\lambda}^3)\|\varphi_1(\lambda,x)\|_{\mathcal{H}}^2=-i\int\limits_{-\infty}^0[Dv_1(\lambda,x)\cdot\overline{v_1(\lambda,x)}-v_1(\lambda,x)\overline{Dv_1(
\lambda,x)}]dx$$
$$+\int\limits_0^\infty\{[[iD^3+q(x)]e_1(\lambda,x)]\cdot\overline{e_1(\lambda,x)}-e_1(\lambda,x)\overline{[[iD^3+q(x)]e_1(\lambda,x)]}\}dx$$
$$=-i|v_1(\lambda,0)|^2-i\{e''_1(\lambda,0)\overline{e_1(\lambda,0)}-e'_1(\lambda,0)\overline{e'_1(\lambda,0)}+e_1(\lambda,0)\overline{e''_1(\lambda,0)}\}=0,$$
and thus $\lambda^3-\overline{\lambda}^3=0$. Since $\lambda-\overline{\lambda}\not=0$, for $\lambda\in\Omega_2$, then $\lambda^2+\lambda\overline{\lambda}+\overline{\lambda}^2=0$, hence it follows that $\lambda=\overline{\lambda}\zeta_2$ or $\lambda=\overline{\lambda}\zeta_3$. For $\lambda=x\eta$ ($x>0$, $\eta\in\mathbb{T}$), we have $\eta^2=\zeta_2$ or $\eta^2=\zeta_3$. Therefore, either $\lambda=-\varkappa\zeta_3$ and a negative eigenvalue corresponds to it, $\lambda^3=-\varkappa^3<0$, or $\lambda=\varkappa\zeta_2$, to which there corresponds a positive $\lambda^3=\varkappa^3>0$ eigenvalue. As a result, we have two series of roots of the function $e_1(\lambda,0)$ in the sector $\Omega_2$,
$$\lambda_n=-\varkappa_n\zeta_3\in\widehat{l}_{\zeta_3},\quad\mu_n=\varkappa_n\zeta_2\in l_{\zeta_2}\quad(\varkappa_n>0,n\in\mathbb{N}).$$
The component $v_1(\lambda,x)$ of the function $\varphi_1(\lambda,x)$ \eqref{eq2.36}, for $\lambda=\lambda_n$, is equal to $v_1(\lambda_n,x)=\theta e'_1(\lambda_n,0)\cdot e^{-i\varkappa_n^3x}$ and does not belong to $L^2(\mathbb{R}_-)$ if $e'_1(\lambda_n,0)\not=0$. Hence, as was noted before, the eigenfunction $\varphi_1(\lambda_n,x)$ belongs to $\mathcal{H}$ only in the case when $\lambda_n$ is a joint zero of the functions $e_1(\lambda,0)$, $e'_1(\lambda,0)$.

Prove that the sequences $\{\lambda_n\}$ and $\{\mu_n\}$ are finite. Assuming the contrary, suppose that, e. g., the set $\{\lambda_n\}$ is countable. From $\{\lambda_n\}$, we single out a converging subsequence $\{\lambda_{n_k}\}$, such that $\lambda_{n_k}\rightarrow0$ ($n_k\rightarrow\infty$). The eigenfunctions $\varphi_1(\lambda_{n_k},x)$ \eqref{eq2.36} are given by $\varphi_1(\lambda_{n_k},x)=(0,e_1(\lambda_n,x))$ and $\varphi_1(\lambda_{n_k},x)\stackrel{w}{\rightarrow}\varphi_1(0,x)$. Using the orthogonality $\varphi_1(\lambda_{n_k},x)\perp\varphi_1(\lambda_{n_s},x)$ ($\lambda_{n_k}\not=\lambda_{n_s}$), we obtain
$$0=\langle\varphi_1(\lambda_{n_k},x),\varphi_1(\lambda_{n_s},x)\rangle=\langle\varphi_1(\lambda_{n_k},x),\varphi_1(0,x)\rangle$$
$$+\langle\varphi_1(\lambda_{n_k},x),[\varphi_1(
\lambda_{n_s},x)-\varphi_1(0,x)].$$
Tending $\lambda_{n_s}\rightarrow0$ in this equality, we have $0=\langle\varphi_1(\lambda_{n_k},x),\varphi_1(0,x)\rangle$ and after passing to the limit $\lambda_{n_k}\rightarrow0$ we find that $\|\varphi_1(0,x)\|=0$ and thus $e_1(0,x)=0$. And since (see \eqref{eq2.8})
$$e_1(0,x)=1-\int\limits_x^\infty\frac{(x-t)^2}2q(t)e_1(0,t)dt,$$
then substituting $e_1(0,x)\equiv0$, we arrive at the contradiction, $0=1$.

Prove that zeros $\lambda\in\Lambda_1$ \eqref{eq2.35} of the function $e_1(\lambda,0)$ are simple. Upon differentiating by $\lambda$ the equality
\begin{equation}
ie'''_1(\lambda,x)+q(x)e_1(\lambda,x)=\lambda^3e_1(\lambda,x)\quad(\lambda\in\Omega_2)\label{eq2.37}
\end{equation}
and supposing that $\partial_\lambda y(\lambda,x)\stackrel{\rm def}{=}dy(\lambda,x)/d\lambda$, we have
$$i\partial_\lambda e'''_1(\lambda,x)+q(x)\partial_\lambda e_1(\lambda,x)=\lambda^3\partial_\lambda e_1(\lambda,x)+3\lambda^2e_1(\lambda,x).$$
Hence, after taking complex conjugate, in view of the fact that $q(x)$ and $\lambda^3$ ($\lambda\in\Lambda_1$ \eqref{eq2.35}) are real, we find
$$-i\overline{\partial_\lambda e''_1(\lambda,x)}+q(x)\overline{\partial_\lambda e_1(\lambda,x)}=\lambda^3\overline{\partial_\lambda e_1(\lambda,x)}+3\overline{\lambda}^2e_1(\lambda,x).$$
Multiplying this equality by $e_1(\lambda,x)$ and \eqref{eq2.37}, by $\overline{\partial_\lambda e_1(\lambda,x)}$, and subtracting the obtained, we have
$$ie'''_1(\lambda,x)\overline{\partial_\lambda e_1(\lambda,x)}+e_1(\lambda,x)\overline{\partial_\lambda e'''_1(\lambda,x)}=-3\overline{\lambda}^2|e_1(\lambda,x)|^2.$$
Integration by $x$ from 0 to $\infty$, due to $e(\lambda,0)=e'(\lambda,0)=0$, gives
$$ie''_1(\lambda,0)\overline{\partial_\lambda e_1(\lambda,0)}=3\overline{\lambda}^2\int\limits_0^\infty |e_1(\lambda,x)|^2dx,$$
and thus $\partial_1e_1(\lambda,0)\not=0$. $\blacksquare$

\begin{picture}(200,200)
\put(0,100){\vector(1,0){200}}
\put(100,0){\vector(0,1){200}}
\put(150,0){\vector(-1,2){100}}
\put(150,200){\vector(-1,-2){100}}
\put(200,67){\line(-3,1){200}}
\put(0,67){\line(3,1){200}}
\put(30,190){$l_{\zeta_2}$}
\put(135,150){$\lambda_n$}
\put(80,170){$\Omega_2$}
\put(102,190){$e_1(\lambda,x)$}
\put(150,180){$\widehat{l}_{\zeta_3}$}
\put(121,150){$\times$}
\put(70,147){$\circ$}
\put(150,96){$\circ$}
\put(154,105){$\mu'_n$}
\put(140,50){$e_2(\lambda,x)$}
\put(140,30){$\widehat{l}_{\zeta_2}$}
\put(122,40){$\times$}
\put(110,35){$\lambda'_n$}
\put(68,40){$\circ$}
\put(80,35){$\mu''_n$}
\put(60,5){$l_{\zeta_3}$}
\put(60,60){$\Omega_3$}
\put(10,85){$e_3(\lambda,x)$}
\put(10,105){$\widehat{l}_{\zeta_1}$}
\put(50,96){$\times$}
\put(38,105){$\lambda''_n$}
\put(190,105){$l_{\zeta_1}$}
\put(170,85){$\Omega_1$}
\put(50,150){$\mu_n$}
\qbezier(123,109)(100,130)(82,106)
\qbezier(100,80)(120,90)(120,107)
\qbezier(80,107)(80,80)(100,77)
\end{picture}

\hspace{20mm} Fig. 6

For the functions $e_2(\lambda,x)$ and $e_3(\lambda,x)$, holomorphic in the sectors $\Omega_1$ and $\Omega_3$ \eqref{eq2.33} correspondingly, an analogue of Theorem \ref{t2.2} is true. Really, using \eqref{eq2.23}, we arrive at the following description of joint zeros of the pairs of functions $e_2(\lambda,0)$, $e'_2(\lambda,0)$ and $e_3(\lambda,0)$, $e'_3(\lambda,0)$. These sets are given by
\begin{equation}
\begin{array}{lll}
\Lambda_2\stackrel{\rm def}{=}\{\lambda'_n,\mu'_n\in\Omega_1:\lambda'_n=\zeta_3\lambda_n=-\varkappa_n\zeta_2\in\widehat{l}_{\zeta_2};\mu'_n=\zeta_3\mu_n=\varkappa_n\zeta_1\in l_{\zeta_1}\};\\
\Lambda_3\stackrel{\rm def}{=}\{\lambda''_n,\mu''_n\in\Omega_3:\lambda''_n=\zeta_2\lambda_n=-\varkappa_n\zeta_1\in\widehat{l}_{\zeta_1};\mu''_n=\zeta_2\mu_n=\varkappa_n\zeta_3\in l_{\zeta_3}\}.
\end{array}\label{eq2.38}
\end{equation}

\begin{remark}\label{r2.4}
So, it is sufficient to know the symmetric set of zeros $\lambda''_n=\varkappa_n\zeta_1$, $\mu'_n=\varkappa_n\zeta_1$ ($\lambda''_n=-\mu'_n$) on one of the straight lines $\{L_{\zeta_p}\}$ \eqref{eq1.8}, e. g., on $L_{\zeta_1}$, and zeros on other lines $L_{\zeta_2}$, $L_{\zeta_3}$ are derived by $\zeta_2$ and $\zeta_3$ rotations; $\zeta_2L_{\zeta_1}=L_{\zeta_2}$, $\zeta_3L_{\zeta_1}=L_{\zeta_3}$.
\end{remark}

\begin{corollary}\label{c2.1}
Operator $\mathcal{L}_q$ \eqref{eq2.2} -- \eqref{eq2.4} can have no more than finite number of eigenfunctions $E_1(\lambda,x)$ and $E_2(\lambda,x)$, $E_3(\lambda,x)$ \eqref{eq2.24} (`connected states', \cite{15,1,4,3,2,5}) enumerated by the numbers $\Lambda_1$ \eqref{eq2.35} and $\Lambda_2$, $\Lambda_3$ \eqref{eq2.38} lying in the sectors $\Omega_2$ and $\Omega_1$, $\Omega_3$ \eqref{eq2.33}. Eigenvalues $\lambda^3$ of the operator $\mathcal{L}_q$ corresponding to $\lambda_n$, $\lambda'_n$, $\lambda''_n$ are negative, and corresponding to $\mu_n$, $\mu'_n$, $\mu''_n$, positive.
\end{corollary}
\vspace{5mm}

{\bf 2.3} Consider the Cauchy problem on the semiaxis $\mathbb{R}_+$ with initial data at zero,
\begin{equation}
\left\{
\begin{array}{lll}
iD^3u(x)+q(x)u(x)=\lambda^3u(x)\quad(x\in\mathbb{R}_+,\lambda\in\mathbb{C});\\
u(0)=0;\,u'(0)=\alpha,\,u''(0)=\beta\quad(\alpha,\beta\in\mathbb{R}).
\end{array}\right.\label{eq2.39}
\end{equation}
Define the sectors $\{\Omega_p^-\}$ obtained from $\{\Omega_p\}$ \eqref{eq2.33} upon the substitution $\lambda\rightarrow-\lambda$,
\begin{equation}
\Omega_p^-\stackrel{\rm def}{=}\{\lambda\in\mathbb{C}:-\lambda\in\Omega_p\}\quad(1\leq p\leq3).\label{eq2.40}
\end{equation}

\begin{lemma}\label{l2.4}
Solution $w(\lambda,x)$ to the Cauchy problem \eqref{eq2.39} is an entire function of $\lambda$, and for all $\lambda\in\overline{(\Omega_2^-)}\backslash\{0\}$ the following inequalities are true:
\begin{equation}
\begin{array}{lll}
(a){\displaystyle\left|\left(w(\lambda,x)-\alpha\frac{s_1(i\lambda x)}{i\lambda}-\beta\frac{s_2(i\lambda x)}{(i\lambda)^2}\right)e^{-i\lambda x\zeta_1}\right|\leq\frac1{|\lambda|^3}\left(|\alpha|+\frac{|\beta|}{|\lambda|}\right)\cdot p(x);}\\
(b){\displaystyle\left|\left(w'(\lambda,x)-\alpha s_0(i\lambda x)-\frac{\beta s_1(i\lambda x)}{i\lambda}\right)e^{-i\lambda x\zeta_1}\right|\leq\frac1{|\lambda|^2}\left(|\alpha|+\frac{|\beta|}{|\lambda|}\right)\cdot p(x);}\\
(c){\displaystyle|(w''(\lambda,x)-\alpha i\lambda s_2(i\lambda x)-\beta s_0(i\lambda x))e^{-i\lambda x\zeta_1}|\leq\frac1{|\lambda|}\left(|\alpha|+\frac{|\beta|}{|\lambda|}\right)p(x);}\
\end{array}\label{eq2.41}
\end{equation}
where $\overline{(\Omega_2^-)}$ is closure of the sector $\Omega_2^-$ \eqref{eq2.40}, $\sigma(x)$ is given by \eqref{eq2.18}, and
\begin{equation}
p(x)\stackrel{\rm def}{=}\sigma(x)\left(1+\frac{x^2}2\sigma(x)\cdot\exp\{x^2\sigma(x)\}\right).\label{eq2.42}
\end{equation}
\end{lemma}

P r o o f. The Cauchy problem \eqref{eq2.39}, due to \eqref{eq1.6}, \eqref{eq1.7}, is equivalent to the integral equation
\begin{equation}
w(\lambda,x)=\alpha\frac{s_1(i\lambda x)}{i\lambda}+\beta\frac{s_2(i\lambda x)}{(i\lambda)^2}+i\int\limits_0^x\frac{s_2(i\lambda(x-t))}{(i\lambda)^2}q(t)w(\lambda,t)dt,\label{eq2.43}
\end{equation}
and, for $w_1(\lambda,x)\stackrel{\rm def}{=}w(\lambda,x)e^{-i\lambda\zeta_1x}$,
\begin{equation}
w_1(\lambda,x)=v_1(\lambda,x)+i\int\limits_0^xe^{i\lambda(t-x)}\frac{s_2(i\lambda(x-t))}{(i\lambda)^2}q(t)w_1(\lambda,t)dt\label{eq2.44}
\end{equation}
where
\begin{equation}
v_1(\lambda,x)=\left(a\frac{s_1(i\lambda x)}{i\lambda}+b\frac{s_2(i\lambda x)}{(i\lambda)^2}\right)e^{-i\lambda\zeta_1x}.\label{eq2.45}
\end{equation}
Equation \eqref{eq2.44}, in terms of the Volterra operator
\begin{equation}
\begin{array}{ccc}
{\displaystyle(T_\lambda f)(x)\stackrel{\rm def}{=}\int\limits_0^xT_1(\lambda,x,t)q(t)f(t)dt}\\
{\displaystyle\left(T_1(\lambda,x,t)\stackrel{\rm def}{=}e^{i\lambda(x-t)}\frac{s_2(i\lambda(x-t))}{(i\lambda)^2},f\in L^2(\mathbb{R}_+)\right),}
\end{array}\label{eq2.46}
\end{equation}
becomes
$$(I-iT_\lambda)w_1(\lambda,\lambda)=v_1(\lambda,x),$$
and thus
\begin{equation}
w_1(\lambda,x)=\sum\limits_0^\infty i^nT_\lambda^nV_1(\lambda,x).\label{eq2.47}
\end{equation}
The operators $T_\lambda^n$ ($n\in\mathbb{N}$) are also Volterra,
$$(T_\lambda^nf)(x)=\int\limits_0^xT_n(\lambda,x,t)q(t)f(t)dt\quad(f\in L^2(\mathbb{R}_+),$$
and for the kernels $T_n(\lambda,x,t)$ the following recurrent relations are true:
\begin{equation}
T_{n+1}(\lambda,x,t)=\int\limits_t^xT_1(\lambda,x,s)q(s)T_n(\lambda,s,t)ds\quad(n\in\mathbb{N}).\label{eq2.48}
\end{equation}
The kernel $T_1(\lambda,x,t)$ \eqref{eq2.46} equals
\begin{equation}
T_1(\lambda,x,t)=\frac1{3(i\lambda)^2}\{1+\zeta_2e^{i\lambda(\zeta_2-\zeta_1)(x-t)}+\zeta_3e^{i\lambda(\zeta_3-\zeta_1)(x-t)}\},\label{eq2.49}
\end{equation}
then, according to \eqref{eq2.31}, modulo of each exponent in \eqref{eq2.49} is less than one, for $\lambda\in\Omega_2^-$, and these exponents tend to zero when $\lambda\rightarrow\infty$ ($\lambda\in\Omega_2^-$). Using the equalities
$$\frac{e^{i\lambda(\zeta_2-\zeta_1)x}}{\lambda^2}=\frac{e^{\sqrt3\lambda\zeta_3x}}{\lambda^2}=3\zeta_3^2\int\limits_0^x(x-\tau)e^{\sqrt3\lambda\zeta_3\tau}d\tau+
\frac1{\lambda^2}
+\frac{\sqrt3\zeta_3x}\lambda;$$
$$\frac{e^{i\lambda(\zeta_3-\zeta_1)x}}{\lambda^2}=\frac{e^{-\sqrt3\lambda\zeta_2x}}{\lambda^2}=3\zeta_2^2\int\limits_0^x(x-\tau)e^{-\sqrt3\lambda\zeta_2\tau}d\tau+
\frac1{\lambda^2}-\frac{\sqrt3\zeta_2x}\lambda$$
which are true for all $x\in\mathbb{R}$, re-write $T_1(\lambda,x,t)$ \eqref{eq2.49} as
$$T_1(\lambda,x,t)=-\int\limits_0^{x-t}(x-t-\tau)[\zeta_3e^{\sqrt3\lambda\zeta_3\tau}+\zeta_2e^{-\sqrt3\lambda\zeta_2\tau}]d\tau,$$
and thus
\begin{equation}
|T_1(\lambda,x,t)|<(x-t)^2\quad(\forall\lambda\in\Omega_2^-).\label{eq2.50}
\end{equation}
Analogously to Lemma \ref{l2.1}, \eqref{eq2.48}, \eqref{eq2.50} imply the estimations
\begin{equation}
|T_n(\lambda,x,t)|\leq\frac{(x-t)^{2n}}{n^{2n}}\frac{\sigma^{n-1}(x)}{(n-1)!}\quad(\lambda\in\Omega_2^-,n\in\mathbb{N}).\label{eq2.51}
\end{equation}
Relation \eqref{eq2.47} yields
\begin{equation}
w_1(\lambda,x)=v_1(\lambda,x)+\int\limits_0^xR(\lambda,x,t)q(t)v_1(\lambda,t)dt\label{eq2.52}
\end{equation}
where
\begin{equation}
R(\lambda,x,t)\stackrel{\rm def}{=}\sum\limits_1^\infty i^nT_n(\lambda,x,t),\label{eq2.53}
\end{equation}
besides, series \eqref{eq2.53} uniformly converges and, due to \eqref{eq2.51},
\begin{equation}
|R(\lambda,x,t)|\leq(x-t)^2\cdot\exp\{(x-t)^2\sigma(x)\}\quad(\forall\lambda\in\Omega_2^-).\label{eq2.54}
\end{equation}
The function $v_1(\lambda,x)$ \eqref{eq2.45} is holomorphic in $\Omega_2^-$ and continuous in the closure $\overline{(\Omega_2^-)}$, which, in view of uniform convergence of series \eqref{eq2.53} (see \eqref{eq2.52}), gives analyticity of $w_1(\lambda,x)$ in the sector $\Omega_2^-$ and continuity in $\overline{(\Omega_2^-)}$. Therefore $w(\lambda,x)$ ($=w_1(\lambda,x)e^{i\lambda\zeta_1x}$) is holomorphic in $\Omega_2^-$ and continuous in $\overline{(\Omega_2^-)}$.

Proceed to the proof of inequalities \eqref{eq2.41}. Use the estimates
\begin{equation}
\left|e^{-i\lambda\zeta_1x}\frac{s_p(i\lambda x)}{(i\lambda)^p}\right|\leq\frac1{|\lambda|^p}\quad(p=1,2,\lambda\in\Omega_2^-\backslash\{0\},x\in\mathbb{R}_+).\label{eq2.55}
\end{equation}
From \eqref{eq2.52}, \eqref{eq2.54}, we find that
\begin{equation}
\begin{array}{ccc}
{\displaystyle|w_1(\lambda,x)|\leq\left(\frac{|\alpha|}{|\lambda|}+\frac{|\beta|}{|\lambda|^2}\right)\left(1+\int\limits_0^x(x-t)^2\exp\{(x-t)^2\sigma(x)\}\cdot|q(t)|dt\right)}\\
{\displaystyle\leq\frac1{|\lambda|}\left(|\alpha|+\frac{|\beta|}{|\lambda|}\right)(1+x^2\sigma(x)\exp\{x^2\sigma(x)\}.}
\end{array}\label{eq2.56}
\end{equation}
Using equation \eqref{eq2.44} and inequalities \eqref{eq2.55}, \eqref{eq2.56}, we obtain
$$|w_1(\lambda,x)-v_1(\lambda,x)|\leq\frac1{|\lambda|^3}\left(|\alpha|+\frac{|\beta|}{|\lambda|}\right)\int\limits_0^x|q(t)|(1+t^2\sigma(t)\cdot\exp\{t^2\sigma(t)\})dt$$
$$\leq\frac1{|\lambda|^3}\left(|\alpha|+\frac{|\beta|}{|\lambda|}\right)\left(\sigma(x)+\frac{x^2}2\sigma^2(x)\cdot\exp\{x^2\sigma(x)\}\right),$$
which proves (a) \eqref{eq2.41}. To prove inequality (b) \eqref{eq2.41}, differentiate equation \eqref{eq2.43},
$$w'(\lambda,x)=\alpha s_0(i\lambda x)+\frac{\beta s_1(i\lambda x)}{i\lambda}+i\int\limits_0^x\frac{s_1(i\lambda(x-t))}{i\lambda}q(t)w(\lambda,t)dt,$$
and thus
$$\left(w'(\lambda,x)-as_0(i\lambda x)-\beta\frac{s_1(i\lambda x)}{i\lambda}\right)e^{-i\lambda\zeta_1x}$$
$$=i\int\limits_0^xe^{i\lambda(t-x)}\frac{s_1(i\lambda(x-t))}{i\lambda}q(t)w(\lambda,t)dt,$$
hence, upon using \eqref{eq2.55}, \eqref{eq2.56}, follows (b) \eqref{eq2.41}. Relation (c) \eqref{eq2.41} is proved analogously. $\blacksquare$

\begin{remark}\label{r2.5}
In other sectors $\{\Omega_p^-\}$ \eqref{eq2.40}, estimates analogous to \eqref{eq2.41} hold. They follow from \eqref{eq2.41}, upon substitutions $\lambda\rightarrow\lambda\zeta_2$, $\lambda\rightarrow\lambda\zeta_3$, and the equality
\begin{equation}
w(\lambda\zeta_2,x)=w(\lambda,x).\label{eq2.57}
\end{equation}
\end{remark}
\vspace{5mm}

\section{Scattering problem}\label{s3}

{\bf 3.1} Using linear independence of Jost solutions $\{e_p(\lambda,x)\}$ \eqref{eq2.19} (Lemma \ref{l2.2}), re-write the solution $w(\lambda,x)$ to the Cauchy problem \eqref{eq2.39} as
\begin{equation}
w(\lambda,x)=B_1(\lambda)e_1(\lambda,x)+B_2(\lambda)e_2(\lambda,x)+B_3(\lambda)e_3(\lambda,x).\label{eq3.1}
\end{equation}
Equation \eqref{eq2.57} implies that
\begin{equation}
B_p(\lambda\zeta_2)=B_{p'}(\lambda)\quad p'=(p+1)\mod 3\,(1\leq p\leq3).\label{eq3.2}
\end{equation}
Define the operation '$*$' by the formula
\begin{equation}
f^*(\lambda)\stackrel{\rm def}{=}\overline{f(\overline{\lambda})}.\label{eq3.3}
\end{equation}

\begin{remark}\label{r3.1}
The operation '$*$' does not commute with the substitution $\lambda\rightarrow\lambda\zeta_2$ since
$$(f(\lambda\zeta_2))^*=f(\lambda\zeta_3);\quad \left.f^*(\lambda)\right|_{\lambda=\lambda\zeta_2}=f^*(\lambda\zeta_2).$$
\end{remark}

Denote by $W_{p,s}(\lambda,x)$ the Wronskian of functions $e_p(\lambda,x)$, $e_s(\lambda,x)$,
\begin{equation}
\begin{array}{ccc}
W_{p,s}(\lambda,x)\stackrel{\rm def}{=}\{e_p(\lambda,x),e_s(\lambda,x)\}=e_p(\lambda,x)e'_s(\lambda,x)-e'_p(\lambda,x)e_s(\lambda,x)\\
1\leq p,s\leq3).
\end{array}\label{eq3.4}
\end{equation}

\begin{lemma}\label{l3.1}
The Wronskians $\{W_{p,s}(\lambda,x)\}$ \eqref{eq3.4} are given by
\begin{equation}
\begin{array}{ccc}
W_{1,2}(\lambda,x)=\sqrt3\lambda\zeta_3e_2^*(\lambda,x);\quad W_{2,3}(\lambda,x)=\sqrt3\lambda\zeta_1e_1^*(\lambda,x);\\ W_{3,1}(\lambda,x)=\sqrt3\lambda\zeta_2e_3^*(\lambda,x).
\end{array}\label{eq3.5}
\end{equation}
\end{lemma}

P r o o f. Evidently,
$$W'_{p,s}(\lambda,x)=e_p(\lambda,x)e''_s(\lambda,x)-e''_p(\lambda,x)e_s(\lambda,x);$$
$$W''_p(\lambda,x)=e'_p(\lambda,x)e''_s(\lambda,x)-e''_p(\lambda,x)e'_s(\lambda,x).$$
Differentiating again the last equality and taking into account that $\{e_p(\lambda,x)\}$ are solutions to \eqref{eq2.6}, we obtain that $W_{p,s}(\lambda,s)$ satisfies the equation
\begin{equation}
iy'''(x)-q(x)y(x)=-\lambda^3y(x)\label{eq3.6}
\end{equation}
which, upon taking complex conjugate and substitution $\lambda\rightarrow\overline{\lambda}$, coincides with \eqref{eq2.6}. To prove \eqref{eq3.5}, one has to take into account the asymptotic of \eqref{eq2.7}. So, for $W_{1,2}(\lambda,x)$ as $x\rightarrow\infty$, we have $W_{1,2}(\lambda,x)\rightarrow i\lambda(\zeta_2-\zeta_1)e^{i\lambda(\zeta_1+\zeta_2)x}=\sqrt3\lambda\zeta_3e^{-i\lambda\zeta_3}$ (see \eqref{eq2.30}), which means that $W_{1,2}(\lambda,x)=\sqrt3\lambda\zeta_3e_2^*(\lambda,x)$. $\blacksquare$

In order to find the coefficients $\{B_p(\lambda)\}$ in \eqref{eq3.1}, consider the system of equations
$$w(\lambda,x)=\sum\limits_pB_p(\lambda)e_p(\lambda,x);\quad w'(\lambda,x)=\sum\limits_pB_p(\lambda)e'_p(\lambda,x);$$
$$w''(\lambda,x)=\sum\limits_pB_p(\lambda)e''_p(\lambda,x).$$
Determinant $\Delta(\lambda,x)$ \eqref{eq2.26} of this system is $\Delta(\lambda,x)=-3\sqrt3\lambda^3$. Determinant
$$\Delta_1(\lambda,x)\stackrel{\rm def}{=}\det\left[
\begin{array}{ccc}
w(\lambda,x)&e_2(\lambda,x)&e_3(\lambda,x)\\
w'(\lambda,x)&e'_2(\lambda,x)&e'_3(\lambda,x)\\
w''(\lambda,x)&e''_2(\lambda,x)&e''_3(\lambda,x)
\end{array}\right]=w(\lambda,x)W''_{2,3}(\lambda,x)$$
$$-w'(\lambda,x)W'_{2,3}(\lambda,x)+w''(\lambda,x)W_{2,3}(\lambda,x)$$
does not depend on $x$ since $\Delta'_1(\lambda,x)=0$, because $w(\lambda,x)$ is the solution to equation \eqref{eq2.39} and $W_{2,3}(\lambda,x)$, correspondingly, to \eqref{eq3.6}. Using the initial data \eqref{eq2.39} for the function $w(\lambda,x)$ and \eqref{eq3.5}, we obtain
$$\Delta_1(\lambda,x)=-\alpha W'_{2,3}(\lambda,0)+\beta W_{2,3}(\lambda,0)$$
$$=\sqrt3\lambda\zeta_1(-\alpha e_1^{*'}(\lambda,0)+\beta e_1^*(\lambda,0)),$$
whence we find $B_1(\lambda)=\Delta_1(\lambda,x)/\Delta(\lambda,x)$.

\begin{lemma}\label{l3.2}
Coefficients $\{B_p(\lambda)\}$ from \eqref{eq3.1} are
\begin{equation}
\begin{array}{ccc}
{\displaystyle B_1(\lambda)=\frac{\zeta_1}{3\lambda^2}(\alpha e_1^{*'}(\lambda,0)-\beta e_1^*(\lambda,0));\quad B_2(\lambda)=\frac{\zeta_2}{3\lambda^2}(\alpha e_3^{*'}(\lambda,0)-\beta e_3^*(\lambda,0));}\\
{\displaystyle B_3(\lambda)=\frac{\zeta_3}{3\lambda^2}(\alpha e_2^{*'}(\lambda,0)-\beta e_2^*(\lambda,0))}
\end{array}\label{eq3.7}
\end{equation}
where $\alpha$ and $\beta$ are from boundary conditions \eqref{eq2.39}.
\end{lemma}

\begin{remark}\label{r3.2}
Equations \eqref{eq3.5}, \eqref{eq3.7} imply
\begin{equation}
\begin{array}{ccc}
{\displaystyle w(\lambda,x)=\frac\alpha{3\sqrt3\lambda^3}\det\left[
\begin{array}{ccc}
e_1(\lambda,x)&e_2(\lambda,x)&e_3(\lambda,x)\\
e_1(\lambda,0)&e_2(\lambda,0)&e_3(\lambda,0)\\
e''_1(\lambda,0)&e''_2(\lambda,0)&e''_3(\lambda,0)
\end{array}\right]}\\
{\displaystyle+\frac\beta{3\sqrt3\lambda^3}\det\left[
\begin{array}{ccc}
e_1(\lambda,x)&e_2(\lambda,x)&e_3(\lambda,x)\\
e_1(\lambda,0)&e_2(\lambda,0)&e_3(\lambda,0)\\
e'_1(\lambda,0)&e'_2(\lambda,0)&e'_3(\lambda,0)
\end{array}\right].}
\end{array}\label{eq3.8}
\end{equation}
\end{remark}
\vspace{5mm}

{\bf 3.2} In $\mathcal{H}$ \eqref{eq2.1}, consider the function
\begin{equation}
y(\lambda,x)=(C(\lambda)e^{i\lambda^3x},w(\lambda,x)),\label{eq3.9}
\end{equation}
which is the solution to system \eqref{eq2.5}, and $w(\lambda,x)$ is given by \eqref{eq3.1}. Relation \eqref{eq2.57} and invariancy of system \eqref{eq2.5} with regard to $\lambda\rightarrow\lambda\zeta_2$ imply that $C(\lambda\zeta_2)=C(\lambda)$. If $B_p(\lambda)=0$ ($\forall p$), then $w(\lambda,x)=0$ ($\alpha=\beta=0$); excluding this trivial case, suppose that, e. g., $B_1(\lambda)\not=0$, then \eqref{eq3.9} yields
\begin{equation}
y(\lambda,x)B_1^{-1}(\lambda)=(c_1(\lambda)e^{i\lambda^3x},e_1(\lambda,x)+s_2(\lambda)e_2(\lambda,x)+s_3(\lambda)e_3(\lambda,x)),\label{eq3.10}
\end{equation}
where
\begin{equation}
c_1(\lambda)=C(\lambda)/B_1(\lambda),\quad s_2(\lambda)=B_2(\lambda)/B_1(\lambda),\quad s_3(\lambda)=B_3(\lambda)/B_1(\lambda),\label{eq3.11}
\end{equation}
and, due to \eqref{eq3.7},
\begin{equation}
s_2(\lambda)=\zeta_2\frac{\alpha e^{*'}_3(\lambda,0)-\beta e^*_3(\lambda,0)}{\alpha e^{*'}_1(\lambda,0)-\beta e^*_1(\lambda,0)};\quad s_3(\lambda)=\zeta_3\frac{\alpha e^{*'}_2(\lambda,0)-\beta e_2^*(\lambda,0)}{\alpha e^{*'}_1(\lambda,0)-\beta e_1^*(\lambda,0)}.\label{eq3.12}
\end{equation}

\begin{remark}\label{r3.3}
Upon using \eqref{eq2.7}, equation \eqref{eq3.10} implies
\begin{equation}
y(\lambda,x)B_1^{-1}(\lambda)\rightarrow e^{i\lambda\zeta_1x}+s_2(\lambda)e^{i\lambda\zeta_2x}+s_3(\lambda)e^{i\lambda\zeta_3x}\quad(x\rightarrow\infty)\label{eq3.13}
\end{equation}
($\forall\lambda\in\mathbb{D}_{a/3}\backslash\{0\}$), therefore $s_2(\lambda)$ and $s_3(\lambda)$ are the {\bf scattering coefficients} of an 'incident wave' $e^{i\lambda\zeta_1x}$. The function $c_1(\lambda)$ is said to be the {\bf matching coefficient} of the wave $e^{i\lambda\zeta_1x}$. Scattering and matching coefficients of incident waves $e^{i\lambda\zeta_2x}$ and $e^{i\lambda\zeta_3x}$ are obtained from \eqref{eq3.10}, \eqref{eq3.13}, after the substitutions $\lambda\rightarrow\lambda\zeta_2$ and $\lambda\rightarrow\lambda\zeta_3$.
\end{remark}

Taking into account
\begin{equation}
e_p^*(\lambda\zeta_2,x)=e_{p'}^*(\lambda,x)\quad p'=(p+2)\mod 3\,(1\leq p\leq3),\label{eq3.14}
\end{equation}
relation \eqref{eq3.12} yields that
\begin{equation}
s_2(\lambda\zeta_2)=\zeta_2\frac{\alpha e^{*'}_2(\lambda,0)-\beta e_2^*(\lambda,0)}{\alpha e^{*'}_3(\lambda,0)-\beta e_3^*(\lambda,0)};\quad s_2(\lambda\zeta_3)=\zeta_2\frac{\alpha e_1^{*'}(\lambda,0)-\beta e_1^*(\lambda,0)}{\alpha e^{*'}_2(\lambda,0)-\beta e_2^*(\lambda,0)}.\label{eq3.15}
\end{equation}

\begin{lemma}\label{l3.3}
Functions $s_p(\lambda)$ ($p=1$, $2$) \eqref{eq3.11} (\eqref{eq3.12}) satisfy the equalities
\begin{equation}
\begin{array}{lll}
(i)\,s_2(\lambda)s_2(\lambda\zeta_2)s_2(\lambda\zeta_3)=1;\\
(ii)\,s_2(\lambda\zeta_3)s_3(\lambda)=1.
\end{array}\label{eq3.16}
\end{equation}
\end{lemma}

\begin{remark}\label{r3.4}
The functions $s_2(\lambda)$, $s_3(\lambda)$ \eqref{eq3.12} depend on the numbers $\alpha$, $\beta\in\mathbb{R}$, or (which is the same) on the choice of a solution $w(\lambda,x)$.

Equalities (i), (ii) \eqref{eq3.16} follow from the invariance of system \eqref{eq2.5} with regard to the transformation $\lambda\rightarrow\lambda\zeta_2$.
\end{remark}
\vspace{5mm}

{\bf 3.3} Boundary conditions \eqref{eq2.4} for the function $y(\lambda,x)B_1^{-1}(\lambda)$ \eqref{eq3.10} give
\begin{equation}
\left\{
\begin{array}{lll}
s_2(\lambda)e_2(\lambda,0)+s_3(\lambda)e_3(\lambda,0)+e_1(\lambda,0)=0;\\
s_2(\lambda)e'_2(\lambda,0)+S_3(\lambda)e'_3(\lambda,0)+e'_1(\lambda,0)=\theta c_1(\lambda).
\end{array}\right.\label{eq3.17}
\end{equation}
Multiply the first equality by $e'_p(\lambda,x)$, and the second by $e_p(\lambda,0)$, and subtract ($p=1$, 2, 3), then we obtain\
$$\left\{
\begin{array}{lll}
W_{2,1}(\lambda,0)s_2(\lambda)+W_{3,1}(\lambda,0)s_3(\lambda)=-\theta c_1(\lambda)e_1(\lambda,0);\\
W_{3,2}(\lambda,0)s_3(\lambda)+W_{1,2}(\lambda,0)=-\theta c_1(\lambda)e_2(\lambda,0);\\
W_{2,3}(\lambda,0)s_2(\lambda)+W_{1,3}(\lambda,0)=-\theta c_1(\lambda)e_3(\lambda,0);
\end{array}\right.$$
and, using \eqref{eq3.5}, we have
\begin{equation}
\left\{
\begin{array}{lll}
-\sqrt3\lambda\zeta_3e_2^*(\lambda,0)s_2(\lambda)+\sqrt3\lambda\zeta_2e_3^*(\lambda,0)s_3(\lambda)=-\theta c_1(\lambda)e_1(\lambda,0);\\
-\sqrt3\lambda\zeta_1e_1^*(\lambda,0)s_3(\lambda)+\sqrt3\lambda\zeta_3e_2^*(\lambda,0)=-\theta c_1(\lambda)e_2(\lambda,0);\\
\sqrt3\lambda\zeta_1e_1^*(\lambda,0)S_2(\lambda)-\sqrt3\lambda\zeta_2e_3^*(\lambda,0)=-\theta c_1(\lambda)e_3(\lambda,0).
\end{array}\right.\label{eq3.18}
\end{equation}
The last two equalities in \eqref{eq3.18} yield
\begin{equation}
\begin{array}{ccc}
{\displaystyle s_2(\lambda)=\frac{\sqrt3\lambda\zeta_2e_3^*(\lambda,0)-\theta c_1(\lambda)e_3(\lambda,0)}{\sqrt3\lambda\zeta_1e_1^*(\lambda,0)};}\\
{\displaystyle s_3(\lambda)=\frac{\sqrt3\lambda\zeta_3e_2^*(\lambda,0)+\theta c_1(\lambda)e_2(\lambda,0)}{\sqrt3\lambda\zeta_1e_1^*(\lambda,0)}.}
\end{array}\label{eq3.19}
\end{equation}

\begin{lemma}\label{l3.4}
Representations \eqref{eq3.12} and \eqref{eq3.19} of the function $s_2(\lambda)$ ($s_3(\lambda)$) coincide then and only then, when $\alpha=\theta C$, where $C$ ($=C(\lambda)$) from \eqref{eq3.9} does not depend on $\lambda$, the number $\alpha$ is from boundary conditions \eqref{eq2.39} of the function $w(\lambda,x)$ and $\theta\in\mathbb{T}$.
\end{lemma}

P r o o f. Note that equality $\alpha=\theta C$ coincides with the second boundary condition \eqref{eq2.4} for the function $y(\lambda,x)$ \eqref{eq3.9} ($C(\lambda)=C$ - const). Upon equating \eqref{eq3.12} and \eqref{eq3.19} for $s_2(\lambda)$, we obtain
\begin{equation}
\zeta_2\frac{\alpha e_3^{*'}(\lambda,0)-\beta e_3^*(\lambda,0)}{\alpha e_1^{*'}(\lambda,0)-\beta e_1^*(\lambda,0)}=\frac{\sqrt3\lambda\zeta_2e_3^*(\lambda,0)-\theta c_1(\lambda)e_3(\lambda,0)}{\sqrt3\lambda\zeta_1e_1^*(\lambda,0)}\label{eq3.20}
\end{equation}
or
$$\alpha\sqrt3\lambda\zeta_2(W_{1,3}(\lambda,0))^*=-\theta c_1(\lambda)e_3(\lambda,0)(\alpha e_1^{*'}(\lambda,0)-\beta e_1^*(\lambda,0)),$$
and using \eqref{eq3.5} we have $3\lambda^2\alpha=\theta c_1(\lambda)(\alpha e_1^{*'}(\lambda,0)-\beta e_1^*(\lambda,0))$, hence, due to \eqref{eq3.7} and $c_1(\lambda)=C/B_1(\lambda)$
\eqref{eq3.11}, we find that $\alpha=\theta C$. $\blacksquare$

\begin{remark}\label{r3.5}
Relation \eqref{eq3.20} implies that the set of zeros of $e'_1(\lambda,0)$ belongs to the set of zeros of the function $e_1(\lambda,0)$, moreover, zeros of $e'_1(\lambda,0)$ are simple.
\end{remark}

Consider the vector function $f(\lambda)=\col[e_1(\lambda,0),e_2(\lambda,0),e_3(\lambda,0)]$, then system \eqref{eq3.18} becomes
\begin{equation}
\sqrt3\lambda T(\lambda)f^{[*]}(\lambda)=-\theta c_1(\lambda)f(\lambda)\label{eq3.21}
\end{equation}
where
$$T(\lambda)\stackrel{\rm def}{=}\left[
\begin{array}{ccc}
0&-\zeta_3s_2(\lambda)&\zeta_2s_3(\lambda)\\
-s_3(\lambda)&\zeta_3&0\\
s_2(\lambda)&0&-\zeta_2
\end{array}\right]$$
and $f^{[*]}(\lambda)$ is the vector column obtained from $f(\lambda)$ after application of the operation `$*$' \eqref{eq3.3} to its every element. Relation \eqref{eq3.21} implies that
$$f^{[*]}(\lambda)=-\frac{\sqrt3\lambda}{\overline{\theta}C_1^*(\lambda)}T^{[*]}(\lambda)f(\lambda),$$
besides, $T^{[*]}(\lambda)$ is formed from $T(\lambda)$ after application of the operation `$*$' \eqref{eq3.3} to its every element. Substituting this expression into \eqref{eq3.21}, we have
\begin{equation}
T(\lambda)T^{[*]}(\lambda)f(\lambda)=\frac{c_1^*(\lambda)c_1(\lambda)}{3\lambda^2}f(\lambda).\label{eq3.22}
\end{equation}
Taking into account
$$T(\lambda)T^{[*]}(\lambda)$$
$$=\left[
\begin{array}{ccc}
\zeta_3s_2(\lambda)s_3^*(\lambda)+\zeta_2s_3(\lambda)s_2^*(\lambda)&-s_2(\lambda)&-s_3(\lambda)\\
-\zeta_3s_3^*(\lambda)&\zeta_2s_3(\lambda)s_2^*(\lambda)+1&-\zeta_3s_3(\lambda)s_3^*(\lambda)\\
-\zeta_2s_2^*(\lambda)&-\zeta_2s_2(\lambda)s_2^*(\lambda)&\zeta_3s_2(\lambda)s_3^*(\lambda)+1
\end{array}\right],$$
equation \eqref{eq3.22} for the first coordinate of the vector $f(\lambda)$ implies the equality
$$(\zeta_3s_2(\lambda)s_3^*(\lambda)+\zeta_2s_3(\lambda)s_2^*(\lambda))e_1(\lambda,0)-s_2(\lambda)e_2(\lambda,0)-s_3(\lambda)e_3(\lambda,0)$$
$$=\frac{c_1^*(\lambda)c_1(\lambda)}{3
\lambda^2}e_1(\lambda,0).$$
Hence, upon using the first equation in \eqref{eq3.17}, it follows that
\begin{equation}
\zeta_3s_2(\lambda)s_3^*(\lambda)+\zeta_2s_3(\lambda)s_2^*(\lambda)+1=\frac{c_1(\lambda)c_1^*(\lambda)}{3\lambda^2}.\label{eq3.23}
\end{equation}
Equating the other components of the vectors from equality \eqref{eq3.22} again gives relation \eqref{eq3.23}.

\begin{lemma}\label{l3.5}
For the scattering coefficients $s_2(\lambda)$ and $s_3(\lambda)$, for all $\lambda$, identity \eqref{eq3.23} holds, where $c_1(\lambda)=C/B_1(\lambda)$ \eqref{eq3.11}.
\end{lemma}

\begin{remark}\label{r3.6}
One has to consider relation \eqref{eq3.23} as unitarity property of the studied scattering problem. In this case, it is an analogue of the well-known unitarity property of a scattering matrix \cite{1} -- \cite{5}, \cite{19}. Characteristically, unitarity property of a scattering problem is a corollary of self-adjoint boundary conditions \eqref{eq2.4} of the main operator $\mathcal{L}$ \eqref{eq2.2}.
\end{remark}
\vspace{5mm}

{\bf 3.4} Consider the function
\begin{equation}
\omega_p(\lambda,x)\stackrel{\rm def}{=}w(\lambda,x)/B_p^{-1}(\lambda)\quad(1\leq p\leq3),\label{eq3.24}
\end{equation}
where $w(\lambda,x)$ is given by \eqref{eq3.1} and $B_p(\lambda)$, correspondingly, by \eqref{eq3.7}. It is evident that
\begin{equation}
\omega_p(\lambda\zeta_2,x)=\omega_{p'}(\lambda,x),\quad p'=(p+1)\mod 3\,(1\leq p\leq3).\label{eq3.25}
\end{equation}
And let $\omega_{p,s}(\lambda,x)$ be the Wronskian of the functions $\omega_p(\lambda,x)$ and $e_3(\lambda,x)$,
\begin{equation}
\begin{array}{ccc}
\omega_{p,s}(\lambda,x)\stackrel{\rm def}{=}\{\omega_p(\lambda,x),e_s(\lambda,x)\}\\
=\omega_p(\lambda,x)e'_s(\lambda,x)-\omega'_p(\lambda,x)e_s(\lambda,x)\quad(1\leq p,s\leq3).
\end{array}\label{eq3.26}
\end{equation}
Multiply the first equation of the system
\begin{equation}
\left\{
\begin{array}{lll}
\omega_1(\lambda,x)=e_1(\lambda,x)+s_2(\lambda)e_2(\lambda,x)+s_3(\lambda)e_3(\lambda,x);\\
\omega'_1(\lambda,x)=e'_1(\lambda,x)+s_2(\lambda)e'_2(\lambda,x)+s_3(\lambda)e'_3(\lambda,x);
\end{array}\right.\label{eq3.27}
\end{equation}
by $e'_3(\lambda,x)$ and the second, by $e_3(\lambda,x)$, then, upon subtraction, we obtain
$$\omega_{1,3}(\lambda,x)=W_{1,3}(\lambda,x)+s_2(\lambda)W_{2,3}(\lambda,x).$$
Using \eqref{eq3.5}, we arrive at the statement.

\begin{lemma}\label{l3.6}
For all $\lambda\in\mathbb{D}_{a/3}\backslash\{0\}$, the following equalities hold:
\begin{equation}
\begin{array}{lll}
{\displaystyle ({\rm i})\,s_2(\lambda)e_1^*(\lambda,x)=\zeta_2e_3^*(\lambda,x)+\frac{\zeta_1}{\sqrt3\lambda}\omega_{1,3}(\lambda,x);}\\
{\displaystyle ({\rm ii})\,s_2(\lambda\zeta_2)e_3^*(\lambda,x)=\zeta_2e_2^*(\lambda,x)+\frac{\zeta_3}{\sqrt3\lambda}\omega_{2,1}(\lambda,x);}\\
{\displaystyle ({\rm iii})\,s_2(\lambda\zeta_3)e_2^*(\lambda,x)=\zeta_2e_1^*(\lambda,x)+\frac{\zeta_2}{\sqrt3\lambda}\omega_{3,2}(\lambda,x);}
\end{array}\label{eq3.28}
\end{equation}
where $\{\omega_{p,s}(\lambda,x)\}$ are given by \eqref{eq3.26}.
\end{lemma}

Analogously, upon multiplication of the first equation in \eqref{eq3.27} by $e'_2(\lambda,x)$ and the second, by $e_2(\lambda,x)$, and subtracting, we have
$$\omega_{1,2}(\lambda,x)=W_{1,2}(\lambda,x)+s_3(\lambda)W_{3,2}(\lambda,x),$$
hence, we obtain the following lemma.

\begin{lemma}\label{l3.7}
For all $\lambda\in\mathbb{D}_{a/3}\backslash\{0\}$, the following relations are true:
\begin{equation}
\begin{array}{lll}
{\displaystyle ({\rm i})\,s_3(\lambda)e_1^*(\lambda,x)=\zeta_3e_2^*(\lambda,x)-\frac{\zeta_1}{\sqrt3\lambda}\omega_{1,2}(\lambda,x);}\\
{\displaystyle ({\rm ii})\,s_3(\lambda\zeta_2)e_3^*(\lambda,x)=\zeta_3e_1^*(\lambda,x)-\frac{\zeta_3}{\sqrt3\lambda}\omega_{2,3}(\lambda,x);}\\
{\displaystyle ({\rm iii})\,s_3(\lambda\zeta_3)e_2^*(\lambda,x)=\zeta_3e_3^*(\lambda,x)-\frac{\zeta_2}{\sqrt3\lambda}\omega_{3,1}(\lambda,x);}
\end{array}\label{eq3.29}
\end{equation}
where $\{\omega_{p,s}(\lambda,x)\}$ are given by \eqref{eq3.26}.
\end{lemma}

\begin{remark}\label{r3.7}
Relations \eqref{eq3.28} and \eqref{eq3.29} are equivalent. Moreover, all the equalities in \eqref{eq3.28}, \eqref{eq3.29} follow from one of the equalities in \eqref{eq3.28} (or \eqref{eq3.29}). So, relations $({\rm ii})$, $({\rm iii})$ \eqref{eq3.28} ($({\rm ii})$, $({\rm iii})$ \eqref{eq3.29}) follow from equality $({\rm i})$ \eqref{eq3.28} ($({\rm i})$ \eqref{eq3.29}) upon the substitutions $\lambda\rightarrow\lambda\zeta_2$, $\lambda\rightarrow\lambda\zeta_3$. Show that $({\rm i})$ \eqref{eq3.28} coincides with $({\rm ii})$ \eqref{eq3.29}. To do this, taking into account \eqref{eq3.11}, re-write $({\rm i})$ \eqref{eq3.28} as
$$B_2(\lambda)e_1^*(\lambda,x)-\zeta_2B_1(\lambda)e_3^*(\lambda,x)=\frac{\zeta_1}{\sqrt3\lambda}[w(\lambda,x)e'_3(\lambda,x)-w'(\lambda,x)e_3(\lambda,x)].$$
Equality $({\rm ii})$ \eqref{eq3.29}, due to $s_3(\lambda\zeta_2)=B_1(\lambda)/B_2(\lambda)$ (see \eqref{eq3.2}, \eqref{eq3.11}), becomes
$$B_1(\lambda)e_3^*(\lambda,x)-\zeta_3B_2(\lambda)e_1^*(\lambda,x)=-\frac{\zeta_2}{\sqrt3\lambda}[w(\lambda,x)e'_3(\lambda,x)-w'(\lambda,x)e_3(\lambda,x)],$$
hence follows the identity of the equations $({\rm i})$ \eqref{eq3.28} and $({\rm ii})$ \eqref{eq3.29}.
\end{remark}
\vspace{5mm}

{\bf 3.5.} Equalities \eqref{eq3.28} and \eqref{eq3.29} generate problems of a jump for functions holomorphic in adjacent sectors of a half-plane. Equation (i) \eqref{eq3.29} written in terms of $\{\psi_p(\lambda,x)\}$ \eqref{eq2.27} becomes
\begin{equation}
\zeta_3s_3^*(\lambda)e^{i\lambda(\zeta_1-\zeta_2)x}\psi_1(\lambda,x)=\psi_2(\lambda,x)-f_{1,2}(\lambda,x)\quad\left(f_{1,2}\stackrel{\rm def}{=}\frac{\zeta_3}{\sqrt3\lambda}e^{-i\lambda\zeta_2x}\omega_{1,2}^*(\lambda,x)\right).\label{eq3.30}
\end{equation}

\begin{lemma}\label{l3.8}
The functions $\psi_2(\lambda,x)$ and $f_{1,2}(\lambda,x)$ are holomorphic in the sectors $\Omega_1$ and $\Omega_2\cap\Omega_3^-$ in the right half-plane $(-i\mathbb{C}_+)=\{\lambda\in\mathbb{C}:\Re\lambda>0\}$, besides, $f_{1,2}(\lambda,x)$ have simple poles at the points $\lambda=\lambda_n\in\widehat{l}_{\zeta_3}$ ($\lambda_n\in\Lambda_1$ \eqref{eq2.35}). The functions $\psi_2(\lambda,x)$ and $f_{1,2}(\lambda,x)$ tend to $1$ as $\lambda\rightarrow\infty$ (inside the sectors).
\end{lemma}

P r o o f. The functions $\psi_1(\lambda,x)$ and $\psi_2(\lambda,x)$ are holomorphic in the sectors $\Omega_2$ and $\Omega_1$ (Lemma \ref{l2.3}).

\begin{picture}(200,200)
\put(0,100){\vector(1,0){200}}
\put(100,0){\vector(0,1){200}}
\put(150,200){\line(-1,-2){100}}
\qbezier[90](200,67)(100,100)(0,133)
\put(0,67){\line(3,1){200}}
\put(110,160){$\lambda_n$}
\put(130,140){$f_{1,2}(\lambda,x)$}
\put(150,180){$\widehat{l}_{\zeta_3}$}
\put(121,150){$\times$}
\put(130,65){$\psi_2(\lambda,x)$}
\put(180,120){$(-il_{\zeta_2})$}
\put(50,150){$\psi_1(\lambda,x)$}
\qbezier(100,120)(115,125)(122,109)
\qbezier(100,80)(120,90)(120,107)
\end{picture}

\hspace{20mm} Fig. 7

Show that $\omega_{1,2}^*(\lambda,x)$, and thus $f_{1,2}(\lambda,x)$ \eqref{eq3.30} also, is analytical in the sector $\Omega_2\cap\Omega_3^-$. In the representation
$$\omega_{1,2}^*(\lambda,x)=\frac1{B_1^*(\lambda)}\{w^*(\lambda,x),e_2^*(\lambda,x)\}$$
the function $w^*(\lambda,x)$ is entire and $e_2^*(\lambda,x)$, corresponding to $e^{-i\lambda\zeta_3x}$, is holomorphic in $\Omega_3^-$ \eqref{eq2.30}. Equation \eqref{eq3.7} implies analyticity of $B_1^*(\lambda)$ in the sector $\Omega_2$. This means that $f_{1,2}(\lambda)$ is holomorphic in the sector $\Omega_2\cap\Omega_3^-$, excluding simple poles $\lambda_n\in\Lambda_1$ lying on the ray $\widehat{l}_{\zeta_3}$.

Evidently, $\psi_2(\lambda,x)\rightarrow1$ as $\lambda\rightarrow\infty$ ($\lambda\in\Omega_1$) (see Remark \ref{r2.3}). Show that $f_{1,2}(\lambda,x)\rightarrow1$ when $\lambda\rightarrow\infty$ ($\lambda\in\Omega_2\cap\Omega_3^-$). Really,
$$f_{1,2}(\lambda,x)=\left\{\frac{\zeta_1}{\sqrt3\lambda}e^{i\lambda\zeta_3x}\omega_{1,2}(\lambda,x)\right\}^*$$
$$=\left\{\frac{\zeta_2}{\sqrt3\lambda B_1(\lambda)}[w(\lambda,x)e^{-i\lambda\zeta_3x}\cdot e'_2(\lambda,x)e^{-i\lambda\zeta_2x}-w'(\lambda,x)e^{-i\lambda\zeta_3x}\cdot e_2(\lambda,x)e^{-i\lambda\zeta_2x}]\right\}^*.$$
Taking into account \eqref{eq2.34}, \eqref{eq2.41} ($b=0$) and the form of $B_1(\lambda)$ \eqref{eq3.7} ($b=0$), we obtain
$$f_{1,2}(\lambda,x)
=\left\{\frac{\zeta_2\sqrt3\lambda}{\alpha(-i\lambda+o(1/\lambda))}\left[\frac\alpha{3i\lambda}(1+o(1/\lambda))(i\lambda\zeta_2+o(1/\lambda))\right.\right.$$
$$\left.\left.-\frac\alpha3(1+o(1/2))(1+
o(1/\lambda^2))\right]\right\}^*=\left\{\frac{\zeta_2\lambda((\zeta_2-\zeta_1)+o(1/\lambda))}{\sqrt3(-i\lambda+o(1/\lambda))}\right\}^*.$$
And since $\zeta_2-\zeta_1=-i\sqrt3\zeta_3$, this implies that $f_{1,2}(\lambda,x)\rightarrow1$ as $\lambda\rightarrow\infty$ and $\lambda\in\Omega_2\cap\Omega_3^-$. The case of $b\not=0$ is proved analogously. $\blacksquare$

\begin{corollary}\label{c3.1}
In the right half-plane $(-i\mathbb{C}_+)$ on the ray $(-il_{\zeta_2})$ equality \eqref{eq3.30} defines the jump problem (the jump equals $\zeta_3s_3^*(\lambda)e^{i\lambda(\zeta_1-\zeta_2)x}\psi_1(\lambda,x)$) for holomorphic in adjoint sectors $\Omega_1$ and $\Omega_2\cap\Omega_3^-$ functions $\psi_2(\lambda,x)$ and $f_{1,2}(\lambda,x)$ that tend to $1$ as $\lambda\rightarrow\infty$ inside these sectors. The function $f_{1,2}(\lambda,x)$ has the finite number of simple poles $\lambda_n\in\widehat{l}_{\zeta_3}$.
\end{corollary}

Equality (ii) \eqref{eq3.29} gives the problem of a jump,
\begin{equation}
\zeta_3s_3^*(\lambda\zeta_3)e^{i\lambda(\zeta_3-\zeta_1)x}\psi_3(\lambda,x)=\psi_1(\lambda,x)-f_{2,3}(\lambda,x)\quad\left(f_{2,3}(\lambda,x)\stackrel{\rm def}{=}\frac{\zeta_1e^{-i\lambda\zeta_1x}}{\sqrt3\lambda}\omega_{2,3}^*(\lambda,x)\right)\label{eq3.31}
\end{equation}
on the ray $(-il_{\zeta_3})$ in the half-plane $(-i\zeta_2\mathbb{C}_+)$ for the analytical in adjoint sectors $\Omega_2$ and $\Omega_3\cap\Omega_1^-$ functions $\psi_1(\lambda,x)$ and $f_{2,3}(\lambda,x)$ that tend to 1 as $\lambda\rightarrow\infty$ inside these sectors. The function $f_{2,3}(\lambda,x)$ has the finite number of poles $\lambda''_n\in\widehat{l}_{\zeta_1}$.

Analogously, equality (iii) \eqref{eq3.29} results in the jump problem,
\begin{equation}
\zeta_3s_3^*(\lambda\zeta_2)e^{i\lambda(\zeta_2-\zeta_3)x}\psi_2(\lambda,x)=\psi_3(\lambda,x)-f_{3,1}(\lambda,x)\quad\left(f_{3,1}(\lambda,x)\stackrel{\rm def}{=}\frac{\zeta_2e^{-i\lambda\zeta_3x}}{\sqrt3\lambda}\omega_{3,1}^*(\lambda,x)\right)\label{eq3.32}
\end{equation}
on the ray $(-il_{\zeta_1})$ in the half-plane $(-i\zeta_3\mathbb{C}_+)$ for the holomorphic in adjoint sectors $\Omega_3$ and $\Omega_1\cap\Omega_2^-$ functions $\psi_3(\lambda,x)$ and $f_{3,1}(\lambda,x)$ which, for $\lambda\rightarrow\infty$, tend to 1 inside the sectors. The function $f_{3,1}(\lambda,x)$ has simple poles at the points $\lambda'_n\in\widehat{l}_{\zeta_2}$.

Analogous considerations applied to \eqref{eq3.28} lead to three more problems of jumps in half-planes. So, equality (i) \eqref{eq3.28} written as
\begin{equation}
\begin{array}{ccc}
\zeta_2s_2^*(\lambda)e^{i\lambda(\zeta_1-\zeta_3)x}\psi_1(\lambda,x)=\psi_3(\lambda,x)-g_{1,3}(\lambda,x)\\
{\displaystyle\left(g_{1,3}(\lambda,x)\stackrel{\rm def}{=}-\frac{\zeta_2}{\sqrt3\lambda}e^{-i\lambda\zeta_3x}\omega_{1,3}^*(\lambda,x)\right)}
\end{array}\label{eq3.23}
\end{equation}
gives a problem of jump on the ray $(-il_{\zeta_3})$ in the left half-plane $(i\mathbb{C}_+)=\{\lambda\in\mathbb{C}:\Re\lambda<0\}$. Functions $\psi_3(\lambda,x)$ and $g_{1,3}(\lambda,x)$ are holomorphic in the sectors $\Omega_3$ and $\Omega_1\cap\Omega_2^-$, besides, $g_{1,3}(\lambda,x)$ has simple poles at the points $\mu_n\in l_{\zeta_2}$. Moreover, these functions tend to 1 when $\lambda\rightarrow1$ inside these sectors.

Equality (ii) \eqref{eq3.28} gives a jump problem
\begin{equation}
\begin{array}{ccc}
\zeta_2s_2^*(\lambda\zeta_3)e^{i\lambda(\zeta_3-\zeta_2)x}\psi_3(\lambda,x)=\psi_2(\lambda,x)-g_{2,1}(\lambda,x)\\
{\displaystyle\left(g_{2,1}(\lambda,x)\stackrel{\rm def}{=}-\frac{\zeta_3e^{-i\lambda\zeta_2x}}{\sqrt3\lambda}\omega_{2,1}^*(\lambda,x)\right)}
\end{array}\label{eq3.34}
\end{equation}
on the ray $(-il_{\zeta_1})$ in adjoint sectors $\Omega_1$ and $\Omega_3\cap\Omega_2^-$ of the half-plane $(i\zeta_2\mathbb{C}_+)$ for the functions $\psi_2(\lambda,x)$ and $g_{2,1}(\lambda,x)$ holomorphic in these sectors. These functions tend to 1 as $\lambda\rightarrow\infty$ (inside these sectors), moreover, $g_{2,1}(\lambda,x)$ has simple poles at $\mu''_n\in l_{\zeta_3}$.

Similarly, (iii) \eqref{eq3.28} generates a jump problem,
\begin{equation}
\begin{array}{ccc}
\zeta_2s_2^*(\lambda\zeta_2)e^{i\lambda(\zeta_2-\zeta_1)x}\psi_2(\lambda,x)=\psi_1(\lambda,x)-g_{3,2}(\lambda,x)\\
{\displaystyle\left(g_{3,2}(\lambda,x)\stackrel{\rm def}{=}-\frac{\zeta_1}{\sqrt3\lambda}e^{-i\lambda\zeta_1x}\omega_{3,2}^*(\lambda,x)\right)}
\end{array}\label{eq3.35}
\end{equation}
on the ray $(-il_{\zeta_2})$ in adjoint sectors $\Omega_2$ and $\Omega_1\cap\Omega_3^-$ of the half-plane $(i\zeta_3\mathbb{C}_+)$ for the holomorphic in these sectors functions $\psi_1(\lambda,x)$ and $g_{3,2}(\lambda,x)$ ($g_{3,2}(\lambda,x)$ has simple poles at the points $\mu'_n\in l_{\zeta_1}$) which for $\lambda\rightarrow\infty$ tend to 1 (inside the sectors).
\vspace{5mm}

{\bf 3.6} Consider the jump problems \eqref{eq3.31} and \eqref{eq3.35},
\begin{equation}
\left[
\begin{array}{lll}
\zeta_3s_3^*(\lambda\zeta_3)e^{i\lambda(\zeta_3-\zeta_1)x}\psi_3(\lambda,x)=\psi_1(\lambda,x)-f_{2,3}(\lambda,x)&(\lambda\in(-il_{\zeta_3}));\\
\zeta_2s_2^*(\lambda\zeta_2)e^{i\lambda(\zeta_2-\zeta_1)x}\psi_2(\lambda,x)=\psi_1(\lambda,x)-g_{3,2}(\lambda,x)&(\lambda\in(-il_{\zeta_2}))
\end{array}\right.\label{eq3.36}
\end{equation}
where the function $\psi_1(\lambda,x)$ is holomorphic in $\Omega_2$, and $f_{2,3}(\lambda,x)$ and $g_{3,2}(\lambda,x)$, in the sectors $\Omega_3\cap\Omega_1^-$ and $\Omega_1\cap\Omega_3^-$ correspondingly.

\begin{picture}(200,200)
\put(0,100){\vector(1,0){200}}
\put(100,0){\vector(0,1){200}}
\put(200,67){\vector(-3,1){200}}
\put(0,67){\vector(3,1){200}}
\put(2,130){$(-il_{\zeta_3})$}
\put(100,150){$\psi_1(\lambda,x)$}
\qbezier(100,140)(135,145)(141,115)
\qbezier(140,85)(155,105)(143,116)
\qbezier(100,60)(135,70)(144,85)
\qbezier(100,58)(70,60)(57,82)
\qbezier(56,86)(50,100)(58,115)
\qbezier(100,141)(65,150)(57,116)
\put(10,55){$(-il_{\zeta_2})$}
\put(0,105){$f_{2,3}(\lambda,x)$}
\put(180,135){$(-il_{\zeta_2})$}
\put(170,110){$g_{3,2}(\lambda,x)$}
\put(165,53){$(il_{\zeta_3})$}
\put(105,40){$f_{3,1}(\lambda,x)$}
\put(45,40){$g_{2,1}(\lambda,x)$}
\end{picture}

\hspace{20mm} Fig. 8

Equations \eqref{eq3.31} and \eqref{eq3.7}, \eqref{eq3.26} imply that
$$f_{2,3}(\lambda,x)=\frac{\zeta_2\lambda\sqrt3e^{-i\zeta_1x}}{\alpha e'_3(\lambda,0)-\beta e_3(\lambda,0)}\{w^*(\lambda,x),e_3^*(\lambda,x)\}$$
($\lambda\in\Omega_3\cap\Omega_1^-$), and thus the boundary value of this function on the ray $(il_{\zeta_2})$ is
\begin{equation}
f_{2,3}(i\zeta_2t,x)=\frac{i\zeta_3t\sqrt3e^{\zeta_2tx}}{\alpha e'_1(it,0)-\beta e_1(it,0)}\{w^*(it,x),e_2^*(it,x)\}\quad(t\in\mathbb{R}_+).\label{eq3.37}
\end{equation}
Similarly, the boundary value of $g_{3,2}(\lambda,x)$ on the ray $(il_{\zeta_3})$ is
\begin{equation}
g_{3,2}(i\zeta_3t,x)=-\frac{i\zeta_2t\sqrt3e^{\zeta_3tx}}{\alpha e'_1(it,0)-\beta e_1(it,0)}\{w^*(it,x),e_3^*(it,x)\}\quad(t\in\mathbb{R}_+).\label{eq3.38}
\end{equation}
Holomorphic in the sectors $\Omega_1\cap\Omega_2^-$ and $\Omega_3\cap\Omega_2^-$ functions $f_{3,1}(\lambda,x)$ \eqref{eq3.32} and $g_{2,1}(\lambda,x)$ \eqref{eq3.35} on the rays $(il_{\zeta_3})$ and $(il_{\zeta_2})$ become
$$f_{3,1}(i\zeta_3t,x)=\frac{i\zeta_3t\sqrt3e^{\zeta_2tx}}{\alpha e'_1(it,0)-\beta e_1(it,0)}\{w^*(it,x),e_2^*(it,x)\}\quad(t\in\mathbb{R}_+);$$
$$g_{2,1}(i\zeta_2t,x)=-\frac{i\zeta_2t\sqrt3e^{\zeta_3tx}}{\alpha e'_1(it,0)-\beta e_1(it,0)}\{w^*(it,x),e_3^*(it,x)\}\quad(t\in\mathbb{R}_+).$$
Hence it follows that
\begin{equation}
\left.g_{3,2}(\lambda,x)\right|_{\lambda\in(il_{\zeta_3})}=\left.g_{2,1}(\lambda,x)\right|_{\lambda\in(il_{\zeta_2})};\quad\left.f_{2,3}(\lambda,x)\right|_{\lambda\in(il_{\zeta_2}
)}=\left.f_{3,1}(\lambda,x)\right|_{\lambda\in(il_{\zeta_3})}.\label{eq3.39}
\end{equation}

Symmetry $\lambda\rightarrow-\overline{\lambda}$ about the vertical (imaginary) axis $(il_{\zeta_1})$ sets the one-to-one correspondence between the sectors $\Omega_1\cap\Omega_2^-$ and $\Omega_3\cap\Omega_2^-$, and let $f^+(\lambda)\stackrel{\rm def}{=}f(-\overline{\lambda})$. Then equation \eqref{eq3.39} implies that the values of the functions $g_{3,2}(\lambda,x)$ and $g_{2,1}^+(\lambda,x)$ (as for $f_{2,3}(\lambda,x)$ and $f_{3,1}^+(\lambda,x)$ also) on the ray $(il_{\zeta_3})$ (correspondingly, on $(il_{\zeta_2})$) coincide and are continuous, therefore the functions
\begin{equation}
\begin{array}{lll}
G_2(\lambda,x)\stackrel{\rm def}{=}g_{3,2}(\lambda,x)\chi_{\Omega_1\cap\Omega_3^-}+g_{2,1}^+(\lambda,x)\chi_{\Omega_1\cap\Omega_2^-}\quad(\lambda\in\Omega_1);\\
F_3(\lambda,x)\stackrel{\rm def}{=}f_{2,3}(\lambda,x)\chi_{\Omega_3\cap\Omega_1^-}+f_{3,1}^+(\lambda,x)\chi_{\Omega_3\cap\Omega_2^-}\quad(\lambda\in\Omega_3)
\end{array}\label{eq3.40}
\end{equation}
are holomorphic in the sectors $\Omega_1$ and $\Omega_3$ and tend to 1 when $\lambda\rightarrow\infty$ inside the corresponding sector. The equalities
$$f_{3,1}(-it,x)(=f_{3,1}^+(-it,x))=-\frac{it\sqrt3e^{-t\zeta_3x}}{\alpha e'_2(-it,0)-\beta e_2(-it,0)}\{w^*(-it,x),e_1^*(-it,x)\};$$
$$g_{2,1}(-it,x)(=g_{2,1}^+(-it,x))=\frac{it\sqrt3e^{-t\zeta_2x}}{\alpha e'_3(-it,0)-\beta e_3(-it,0)}\{w^*(-it,x),e_1^*(-it,x)\}$$
imply that
\begin{equation}
G_2(-it,x)=d(-it,x)F_3(-it,x)\label{eq3.41}
\end{equation}
where
\begin{equation}
d(\lambda,x)=-\zeta_2e^{\lambda\sqrt3x}s_2^*(\lambda\zeta_3),\label{eq3.42}
\end{equation}
due to \eqref{eq3.2}, \eqref{eq3.7}, \eqref{eq3.11}. As a result, we arrive at the boundary value Riemann problem \cite{24,25} on the contour $\Gamma\stackrel{\rm def}{=}\bigcup\limits_k(-il_{\zeta_k})$ formed by the rays $(-il_{\zeta_k})$ ($1\leq k\leq3$).

\begin{picture}(200,200)
\put(0,100){\vector(1,0){200}}
\put(100,0){\vector(0,1){200}}
\put(95,0){\vector(0,1){50}}
\put(95,50){\line(0,1){50}}
\put(95,100){\vector(-3,1){50}}
\put(50,115){\line(-3,1){50}}
\put(0,135){\vector(3,-1){50}}
\put(50,119){\line(3,-1){50}}
\put(100,103){\vector(3,1){50}}
\put(150,120){\line(3,1){50}}
\put(200,130){\vector(-3,-1){50}}
\put(150,113){\line(-3,-1){40}}
\put(110,100){\vector(0,-1){50}}
\put(110,50){\line(0,-1){50}}
\put(100,100){\vector(-3,1){100}}
\put(100,100){\vector(3,1){100}}
\put(2,140){$(-il_{\zeta_3})$}
\put(100,140){$\psi_1(\lambda,x)$}
\qbezier(141,115)(100,140)(57,116)
\qbezier(100,60)(140,85)(143,115)
\qbezier(100,58)(57,82)(58,115)
\put(10,55){$F_3(\lambda,x)$}
\put(180,115){$(-il_{\zeta_2})$}
\put(125,50){$G_2(\lambda,x)$}
\put(45,20){$(-il_{\zeta_1})$}
\end{picture}

\hspace{20mm} Fig. 9

On the rays $(-il_{\zeta_3})$ and $(-il_{\zeta_2})$, the boundary value problem coincides with the jump problem \eqref{eq3.36} where $f_{2,3}(\lambda,x)=F_3(\lambda,x)$ as $\lambda\in(-il_{\zeta_3})$ and $g_{3,2}(\lambda,x)=G_2(\lambda,x)$ when $\lambda\in(-il_{\zeta_2})$. And on the ray $(-il_{\zeta_1})$, the boundary condition is given by \eqref{eq3.41}, where $d(\lambda,x)$ is given by \eqref{eq3.42}.

\begin{remark}\label{r3.8}
Function $G_2(\lambda,x)$ \eqref{eq3.40} (and $F_3(\lambda,x)$ \eqref{eq3.40} also) has finite number $(2N)$ of poles at the points $\{\varkappa_n\zeta_1\}\in l_{\zeta_1}$ and $\{-\varkappa_n\zeta_2\}\in\widehat{l}_{\zeta_2}$ ($\{-\varkappa_n\zeta_1\}\in\widehat{l}_{\zeta_1}$ and $\{\varkappa_n\zeta_3\}\in l_{\zeta_3}$ correspondingly), $1\leq n\leq N$. Coefficient $D(\lambda,x)$ of this Riemann problem on the rays $(-il_{\zeta_3})$, $(-il_{\zeta_2})$ equals $1$, and on the ray $(-il_{\zeta_1})$ it coincides with $d(\lambda,x)$ \eqref{eq3.42}.
\end{remark}

Solution of the derived boundary value problem is realized in a standard way \cite{24,25}. First, we find the canonical solution \cite{24,25}
\begin{equation}
\chi(\lambda,x)\stackrel{\rm def}{=}\exp\left\{\frac1{2\pi i}\int\limits_\Gamma\frac{\ln D(\tau,x)}{\tau-\lambda}d\tau\right\}=\exp\left\{-\frac1{2\pi}\int\limits_0^\infty\frac{\ln d(-i\tau,x)}{i\tau+\lambda}d\tau\right\}.\label{eq3.43}
\end{equation}
And next we set the piecewise holomorphic in $\mathbb{C}$ function
\begin{equation}
\mathcal{F}(\lambda,x)\stackrel{\rm def}{=}\left\{
\begin{array}{ccc}
\psi_1(\lambda,x)\chi^{-1}(\lambda,x)&(\lambda\in\Omega_2);\\
G_2(\lambda,x)\chi^{-1}(\lambda,x)&(\lambda\in\Omega_1);\\
F_3(\lambda,x)\chi^{-1}(\lambda,x)&(\lambda\in\Omega_3)
\end{array}\right.\label{eq3.44}
\end{equation}
which is analytic in the sectors $\{\Omega_k\}_1^3$ and tends to $1$ when $\lambda\rightarrow\infty$ (inside sectors). The function $\mathcal{F}(\lambda,x)$ is holomorphic on $(-il_{\zeta_1})$ and the boundary value problem reduces to the problem of jump on the broken line $\widetilde{\Gamma}\stackrel{\rm def}{=}(-il_{\zeta_3})\cup(-il_{\zeta_2})$ (besides, $(-il_{\zeta_3})$ is the incoming ray and $(-il_{\zeta_2})$ is the outgoing ray, see Fig. 9). This problem follows from \eqref{eq3.36} upon division of both sides of the equalities by $\chi(\lambda,x)$. Hence it follows that
\begin{equation}
\mathcal{F}(\lambda,x)=1+b(\lambda,x)+\frac1{2\pi i}\int\limits_{\widetilde{\Gamma}}\frac{\varphi(\tau,x)}{\tau-\lambda}d\tau\label{eq3.45}
\end{equation}
where $\varphi(\tau,x)$ is the jump function on each of the rays of the contour $\widetilde{\Gamma}$ and
\begin{equation}
b(\lambda,x)\stackrel{\rm def}{=}\sum\limits_n\frac{r_n(x)}{\lambda-\varkappa_n}+\sum\limits_n\frac{p_n(x)}{\lambda+\zeta_2\varkappa_n}+\sum\limits_n\frac{\widetilde{r}_n(x)}{\lambda+\varkappa_n}+
\sum\limits_n\frac{\widetilde{p}_n(x)}{\lambda-\zeta_3\varkappa_n},\label{eq3.46}
\end{equation}
besides, $r_n(x)$ and $p_n(x)$ ($\widetilde{r}_n(x)$ and $\widetilde{p}_n(x)$) are residues of the function $G_2(\lambda,x)\times$ $X^{-1}(\lambda,x)$ ($F_3(\lambda,x)X^{-1}(\lambda,x)$) at the points $\lambda=\varkappa_n$ and $\lambda=-\zeta_2\varkappa_n$ ($\lambda=-\varkappa_n$ and $\lambda=\zeta_3\varkappa_n$). Equality \eqref{eq3.45} follows from the Cauchy formula
$$\mathcal{F}(\lambda,1)-1=\frac1{2\pi i}\int\limits_\gamma\frac{\mathcal{F}(\tau,x)}{\tau-\lambda}dt$$
where $\gamma$ is the contour depicted on Fig. 10.

\begin{picture}(400,380)
\put(0,200){\vector(1,0){400}}
\put(200,0){\vector(0,1){380}}
\put(185,202){\vector(-3,1){88}}
\put(100,230){\line(-3,1){30}}
\put(66,250){\vector(3,-1){35}}
\put(100,239){\line(3,-1){100}}
\put(200,200){\line(-1,-2){100}}
\put(200,200){\line(1,-2){100}}
\put(200,206){\vector(3,1){100}}
\put(300,240){\line(3,1){44}}
\put(344,242){\vector(-3,-1){42}}
\put(300,226){\line(-3,-1){70}}
\put(200,200){\line(-3,1){200}}
\put(200,200){\line(3,1){200}}
\put(4,280){$(-il_{\zeta_3})$}
\qbezier(344,256)(200,360)(66,250)
\qbezier(344,242)(200,0)(68,242)
\put(360,240){$(-il_{\zeta_2})$}
\put(155,20){$(-il_{\zeta_1})$}
\put(120,20){$l_{\zeta_3}$}
\put(300,20){$\widehat{l}_{\zeta_2}$}
\put(240,196){$\times$}
\put(280,196){$\times$}
\put(120,196){$\times$}
\put(160,196){$\times$}
\put(179,166){$\times$}
\put(210,166){$\times$}
\put(227,130){$\times$}
\put(161,130){$\times$}
\put(103,185){$-\varkappa_N$}
\put(146,185){$-\varkappa_1$}
\put(330,280){$\gamma$}
\put(5,205){$\widehat{l}_{\zeta_1}$}
\put(380,205){$l_{\zeta_1}$}
\put(230,202){\line(1,0){10}}
\qbezier(240,202)(245,210)(250,202)
\put(250,202){\line(1,0){7}}
\put(257,202){$\dots$}
\put(272,202){\vector(1,0){7}}
\qbezier(279,202)(286,210)(291,200)
\qbezier(291,200)(286,190)(279,198)
\put(279,198){\line(-1,0){7}}
\put(257,198){$\dots$}
\put(257,198){\vector(-1,0){7}}
\qbezier(250,198)(245,190)(240,198)
\put(240,198){\line(-1,0){35}}
\put(205,198){\line(1,-2){12}}
\qbezier(216,174)(223,171)(219,167)
\put(219,167){\line(1,-2){5}}
\multiput(222,155)(2,-4){3}{.}
\put(228,148){\vector(1,-2){5}}
\qbezier(233,138)(243,134)(235,130)
\qbezier(235,130)(223,130)(229,137)
\put(229,137){\line(-1,2){5}}
\multiput(222,147)(-2,4){3}{.}
\put(220,155){\vector(-1,2){5}}
\qbezier(214,164)(205,165)(210,174)
\put(210,174){\line(-1,2){10}}
\put(200,194){\line(-1,-2){10}}
\qbezier(190,175)(195,170)(186,165)
\put(186,165){\line(-1,-2){5}}
\multiput(177,153)(-2,-4){3}{.}
\put(175,145){\vector(-1,-2){5}}
\qbezier(171,136)(175,133)(166,130)
\qbezier(166,130)(155,135)(167,138)
\put(167,138){\line(1,2){5}}
\multiput(169,148)(2,4){3}{.}
\put(174,156){\vector(1,2){5}}
\qbezier(180,166)(173,170)(185,175)
\put(185,175){\line(1,2){11}}
\put(195,197){\line(-1,0){25}}
\qbezier(170,197)(165,192)(160,197)
\put(160,197){\line(-1,0){10}}
\put(135,197){$\dots$}
\put(135,197){\vector(-1,0){5}}
\qbezier(130,197)(125,192)(120,200)
\qbezier(120,200)(124,208)(130,202)
\put(130,202){\line(1,0){5}}
\put(135,202){$\dots$}
\put(150,202){\vector(1,0){9}}
\qbezier(159,202)(164,210)(169,202)
\put(169,202){\line(1,0){17}}
\put(240,185){$\varkappa_1$}
\put(280,185){$\varkappa_N$}
\put(220,165){$-\zeta_2\varkappa_1$}
\put(227,140){$-\zeta_2\varkappa_N$}
\put(175,130){$\zeta_3\varkappa_N$}
\put(150,170){$\zeta_3\varkappa_1$}
\end{picture}

\hspace{20mm} Fig. 10

Integral in \eqref{eq3.45} along the ray $(-il_{\zeta_2})$ ($\tau=-i\zeta_2t$, $t\geq0$) of contour $\gamma$ outgoing from $0$ equals
$$\frac1{2\pi i}\int\limits_{(-il_{\zeta_2})}\frac{\varphi(\tau,x)}{\tau-\lambda}d\tau=\frac1{2\pi i}\int\limits_0^\infty\frac{(-i\zeta_2)dt}{-i\zeta_2t-\lambda}\zeta_2s_2^*(-i\zeta_3t)e^{-i\sqrt3tx}\chi^{-1}(-i\zeta_2t,x)\psi_2(-i\zeta_2t,x)$$
$$=\frac{\zeta_2}{2\pi i}\int\limits_0^\infty\frac{dt}{t-i\zeta_3\lambda}s_2^*(-i\zeta_3t)e^{-i\sqrt3tx}\chi(-i\zeta_2t,x)\psi_2(-i\zeta_2t,x).$$
Analogously, along the incoming ray $(-il_{\zeta_3})$,
$$\frac1{2\pi i}\int\limits_{(-il_{\zeta_3})}\frac{\varphi(\tau,x)}{\tau-\lambda}d\tau=\frac{\zeta_3}{2\pi i}\int\limits_\infty^0\frac{dt}{t-i\zeta_2\lambda}s_3^*(-i\zeta_2t)e^{i\sqrt3tx}\chi^{-1}(-i\zeta_3t,x)\psi_3(-i\zeta_3t,x).$$
Upon substituting these expressions into \eqref{eq3.45} and using $\psi_p(\lambda\zeta_2,x)=\psi_{p'}(\lambda,x)$ ($p'=(p+1)(\mod3)$), we arrive at the representation for $\mathcal{F}(\lambda,x)$ \eqref{eq3.44}
\begin{equation}
\begin{array}{ccc}
{\displaystyle\mathcal{F}(\lambda,x)=1+b(\lambda,x)+\frac{\zeta_2}{2\pi i}\int\limits_0^\infty\frac{d\tau}{\tau-i\zeta_3\lambda}c_2(-i\zeta_2\tau,x)\psi_3(-i\tau,x)}\\
{\displaystyle-\frac{\zeta_3}{2\pi i}\int\limits_0^\infty\frac{d\tau}{\tau-i\zeta_2\lambda}c_3(-i\zeta_3\tau,x)\psi_2(-i\tau,x),}
\end{array}\label{eq3.47}
\end{equation}
here the functions $c_2(\lambda,x)$, $c_3(\lambda,x)$ depend only on $\{s_k(\lambda)\}$,
\begin{equation}
c_2(\lambda,x)\stackrel{\rm def}{=}s_2^*(\lambda\zeta_2)\chi^{-1}(\lambda,x)e^{\sqrt3\zeta_3\lambda x};\quad c_3(\lambda,x)\stackrel{\rm def}{=}s_3^*(\lambda\zeta_3)\chi^{-1}(\lambda,x)e^{-\sqrt3\zeta_2\lambda x},\label{eq3.48}
\end{equation}
besides, $\chi(\lambda,x)$ is given by \eqref{eq3.43} (and $d(\lambda,x)$, by \eqref{eq3.42} correspondingly). Assuming $\lambda\in\Omega_2$ in \eqref{eq3.47} and in view of \eqref{eq3.44}, we obtain
\begin{equation}
\begin{array}{ccc}
{\displaystyle\psi_1(\lambda,x)\chi^{-1}(\lambda,x)=1+b(\lambda,x)+\frac{\zeta_2}{2\pi i}\int\limits_0^\infty\frac{d\tau}{\tau-i\zeta_3x}c_2(-i\zeta_2\tau,x)\psi_3(-i\tau,x)}\\
{\displaystyle-\frac{\zeta_3}{2\pi i}\int\limits_0^\infty\frac{d\tau}{\tau-i\zeta_2\lambda}c_3(-i\zeta_3\tau,x)\psi_2(-i\tau,x).}
\end{array}\label{eq3.49}
\end{equation}
Thus, the holomorphic in $\Omega_2$ function $\psi_1(\lambda,x)\chi^{-1}(\lambda,x)$ is expressed via the functions $\psi_2(-it,x)$, $\psi_3(-it,x)$ and $b(\lambda,x)$ \eqref{eq3.46}. Calculating the boundary values in equation \eqref{eq3.49}, when $\lambda\rightarrow-i\zeta_2t$ and $\lambda\rightarrow-i\zeta_3t$, we obtain the system of integral equations for $\psi_2(-it,x)$ and $\psi_3(-it,x)$,
\begin{equation}
\left\{
\begin{array}{cccc}
{\displaystyle\psi_2(-it,x)\chi^{-1}(-it\zeta_2,x)=1+b(-it\zeta_2,x)+\frac{\zeta_2}{2\pi i}\int\limits_0^\infty\hspace{-4.4mm}/\frac{d\tau}{\tau-t}c_2(-i\zeta_2\tau,x)\psi_3(-i\tau,x)}\\
{\displaystyle-\frac{\zeta_3}{2\pi i}\int\limits_0^\infty\frac{d\tau}{\tau-\zeta_3t}c_3(-i\zeta_3\tau,x)\psi_2(-i\tau,x)+\frac{\zeta_2}2c_2(-i\zeta_2t,x)\psi_3(-it,x);}\\
{\displaystyle\psi_3(-it,x)\chi^{-1}(-it\zeta_3,x)=1+b(-it\zeta_3,x)+\frac{\zeta_2}{2\pi i}\int\limits_0^\infty\frac{d\tau}{\tau-\zeta_2t}c_2(-i\zeta_2\tau,x)\psi_3(-i\tau,x)}\\
{\displaystyle-\frac{\zeta_3}{2\pi i}\int\limits_0^\infty\hspace{-4.4mm}/\frac{d\tau}{\tau-t}c_3(-i\zeta_3\tau,x)-\frac{\zeta_3}2c_3(-i\zeta_3t,x)\psi_2(-it,x).}
\end{array}\right.\label{eq3.50}
\end{equation}
Also it is necessary to derive relations for calculating the functions $\{r_n(x)\}$, $\{p_n(x)\}$, $\{\widetilde{r}_n(x)\}$, $\{\widetilde{p}_n(x)\}$ which are included in $b(\lambda,x)$ \eqref{eq3.46}. Relations \eqref{eq3.31} and \eqref{eq3.35} imply
$$f_{2,3}(\lambda\zeta_2,x)=\frac{\zeta_3}{\sqrt3\lambda}e^{-i\lambda\zeta_2x}\frac{B_3^*(\lambda)}{B_2^*(\lambda\zeta_2)}\omega_{3,2}^*(\lambda,x)=-\zeta_3e^{i\lambda(\zeta_1-
\zeta_2)x}\frac{B_3^*(\lambda)}{B_1^*(\lambda)}g_{3,2}(\lambda,x),$$
therefore
$$\frac{f_{2,3}(\lambda\zeta_2,x)}{B_3^*(\lambda)}=-\zeta_3e^{-\sqrt3\zeta_3\lambda x}\frac{g_{3,2}(\lambda,x)}{B_1^*(\lambda)}.$$
Note that $f_{2,3}(\lambda,x)\chi^{-1}(\lambda,x)=\mathcal{F}(\lambda,x)$ for $\lambda\in\Omega_3\cap\Omega_1^-$ and $g_{3,2}(\lambda,x)\chi^{-1}(\lambda,x)=\mathcal{F}(\lambda,x)$ for $\lambda\in\Omega_1\cap\Omega_3^-$. So,
$$\frac{f_{2,3}(\lambda\zeta_2,x)\chi^{-1}(\lambda\zeta_2,x)}{B_3^*(\lambda)}=-\zeta_3e^{-\sqrt3\zeta_2\lambda x}\chi^{-1}(\lambda\zeta_2,x)\chi(\lambda,x)\cdot\frac{g_{3,2}(\lambda,x)\chi^{-1}(\lambda,x)}{B_1^*(\lambda)},$$
using \eqref{eq3.47}, calculate the residue in both sides of this equality at the point $\lambda=\varkappa_m$ ($a>\varkappa_N$), then
\begin{equation}
\begin{array}{ccc}
{\displaystyle\frac1{B^{*'}_3(\varkappa_m)}\left\{1+b(\zeta_2\varkappa_m,x)+\frac{\zeta_2}{2\pi i}\int\limits_0^\infty\frac{d\tau}{\tau-i\varkappa_m}c_2(-i\zeta_2\tau,x)\psi_3(-i\tau,x)\right.}\\
{\displaystyle\left.-\frac{\zeta_3}{2\pi i}\int\limits_0^\infty\frac{d\tau}{\tau-i\zeta_3\varkappa_m}c_3(-i\zeta_3\tau,x)\psi_2(-i\tau,x)\right\}=}\\
{\displaystyle-\zeta_3e^{-\sqrt3\zeta_3\varkappa_mx}\frac{\chi(\varkappa_m,x)}{\chi(
\zeta_2\varkappa_m,x)B_1^*(\varkappa_m)}r_m(x)}
\end{array}\label{eq3.51}
\end{equation}
($1\leq m\leq N$). Similarly, taking into account
$$\frac{g_{3,2}(\lambda\zeta_3,x)}{B_2^*(\lambda)}=-\zeta_2e^{\sqrt3\zeta_2\lambda x}\frac{f_{2,3}(\lambda,x)}{B_1^*(\lambda)},$$
upon calculating residue at the point $\lambda=-\varkappa_m$, we obtain
\begin{equation}
\begin{array}{ccc}
{\displaystyle\frac1{B_2^{*'}(-\varkappa_m)}\left\{1+b(-\zeta_3\varkappa_m,x)+\frac{\zeta_2}{2\pi i}\int\limits_0^\infty\frac{d\tau}{\tau+i\zeta_2\varkappa_m}c_2(-i\zeta_2\tau,x)\psi_3(-i\tau,x)\right.}\\
{\displaystyle\left.-\frac{\zeta_3}{2\pi i}\int\limits_0^\infty\frac{d\tau}{\tau+i\varkappa_m}c_3(-i\zeta_3\tau,x)\psi_2(-i\tau,x)\right\}}\\
{\displaystyle=-\zeta_2e^{-\sqrt3\zeta_2\varkappa_mx}\frac{\chi(-\varkappa_m,x)}{\chi(-\zeta_3
\varkappa_m,x)B_1^*(-\varkappa_m)}\widetilde{r}_m(x)}
\end{array}\label{eq3.52}
\end{equation}
($1\leq m\leq N$). Equation
$$\frac{f_{2,3}(\lambda\zeta_3,x)}{B_2^*(\lambda)}=-\zeta_3e^{-\sqrt3\zeta_1\lambda x}\frac{g_{2,1}(\lambda,x)}{B_3^*(\lambda)},$$
upon calculation of the residue at $\lambda=\zeta_3\varkappa_m$, implies
\begin{equation}
\begin{array}{ccc}
{\displaystyle\frac1{B_2^{*'}(\zeta_3\varkappa_m)}\left\{1+b(\zeta_2\varkappa_m,x)+\frac{\zeta_2}{2\pi i}\int\limits_0^\infty\frac{d\tau}{\tau-i\varkappa_m}c_2(-i\zeta_2\tau,x)\psi_3(-i\tau,x)\right.}\\
{\displaystyle\left.-\frac{\zeta_3}{2\pi i}\int\limits_0^\infty\frac{d\tau}{\tau-i\zeta_3\varkappa_m}c_3(-i\zeta_3\tau,x)\psi_2(-i\tau,x)\right\}}\\
{\displaystyle=-\zeta_3e^{-\sqrt3\zeta_3\varkappa_mx}\frac{\chi(\zeta_3\varkappa_m,x)}{\chi(\zeta_2
\varkappa_m,x)B_1^*(\varkappa_m)}\widetilde{p}_m(x)}
\end{array}\label{eq3.53}
\end{equation}
($1\leq m\leq N$). Comparing \eqref{eq3.51} and \eqref{eq3.53}, we find that $r_m(x)\chi(\varkappa_m,x)=\chi(\zeta_3\varkappa_m,x)\widetilde{p}_m(x)$ ($1\leq m\leq N$), due to $B_2^{*'}(\zeta_3\varkappa_m)=B_3^{'*}(\varkappa_m)$. Similarly, it is proved that $\chi(-\zeta_2\varkappa_m,x)p_m(x)=\chi(-\varkappa_m,x)\widetilde{r}_m(x)$ ($1\leq m\leq N$). Define the functions
\begin{equation}
z_n(\lambda,x)\stackrel{\rm def}{=}\frac{\chi^{-1}(\varkappa_n,x)}{\lambda-\varkappa_n}+\frac{\chi^{-1}(\zeta_3\varkappa_n,x)}{\lambda-\zeta_3\varkappa_n};\quad\widetilde{z}_n(\lambda,x)\stackrel{\rm def}{=}\frac{\chi^{-1}(-\varkappa_n,x)}{\lambda+\varkappa_n}+\frac{\chi^{-1}(-\zeta_2\varkappa_n,x)}{\lambda+\zeta_2\varkappa_n};\label{eq3.54}
\end{equation}
\begin{equation}
R_n(x)\stackrel{\rm def}{=}r_n(x)\chi(\varkappa_n,x);\quad\widetilde{R}_n(x)=r_n(x)\chi(-\varkappa_n,x)\quad(1\leq n\leq N).\label{eq3.55}
\end{equation}

\begin{lemma}\label{l3.9}
Rational (with respect to $\lambda$) function $b(\lambda,x)$ \eqref{eq3.46} equals
\begin{equation}
b(\lambda,x)=\sum\limits_nR_n(x)z_n(\lambda,x)+\sum\limits_n\widetilde{R}_n(x)\widetilde{z}_n(\lambda,x),\label{eq3.56}
\end{equation}
here $\{R_m(x)\}$, $\{\widetilde{R}_m(x)\}$ are given by \eqref{eq3.55}, functions $\{z_n(\lambda,x)\}$, $\{\widetilde{z}_n(\lambda,x)\}$, correspondingly, by \eqref{eq3.54}, and $\chi(\lambda,x)$ is given by the formula \eqref{eq3.43}.
\end{lemma}

{\bf Conclusion.}{\it Relations \eqref{eq3.50} -- \eqref{eq3.52} ($b(\lambda,x)$ is given by \eqref{eq3.56}) form a {\bf closed system of linear singular equations} relative to $\psi_2(-it,x)$, $\psi_3(-it,x)$ and $\{R_n(x)\}$, $\{\widetilde{R}_n(x)\}$. Thus, solution to the inverse problem is equivalent to the solvability of the equation system \eqref{eq3.50} -- \eqref{eq3.52}.}

\begin{remark}\label{r3.9}
For the third-order operator, we have $2(N+1)$ equations \eqref{eq3.50} -- \eqref{eq3.52}, unlike $N+1$ equations \cite{1,2,6} for the Sturm -- Liouville operator.

Independent parameters that are included in the system \eqref{eq3.50} -- \eqref{eq3.52} are the following: (a) the scattering function $s_2(\lambda)$ ($s_3(\lambda)$ is expressed via $s_2(\lambda)$ by formula (ii) \eqref{eq3.16}); (b) finite set of positive numbers $\{\varkappa_n\}_1^N$, $0<\varkappa_1<...\varkappa_N$ ($0\leq N<\infty$); (c) $2N$ constants $b_n\stackrel{\rm def}{=}B'_3(\varkappa_n)/B_1(\varkappa_n)$, $\widetilde{b}_n\stackrel{\rm def}{=}B'_2(-\varkappa_m)/B_1(-\varkappa_n)$ ($1\leq n\leq N$).
\end{remark}

\section{Inverse problem}\label{s4}

{\bf 4.1} Proceed to solution of the inverse problem. Consider the totality
\begin{equation}
(s_2(\lambda),\{\varkappa_n\}_1^N,\{b_n\}_1^N,\{\widetilde{b}_n\}_1^N)\label{eq4.1}
\end{equation}
where the function $s_2(\lambda)$ is defined on the bundle of straight lines $L$ \eqref{eq1.29}, holomorphic inside circle $\mathbb{D}_r=\{\lambda\in\mathbb{C}:|\lambda|<r\}$ ($0<r<\infty$) and satisfies the condition (i) \eqref{eq3.26}; set of the positive numbers $\{\varkappa_n\}_1^N$ is enumerated in the ascending order, $0<\varkappa_1<...<\varkappa_N$ ($0\leq N<\infty$), and $\{b_n\}_1^N$, $\{\widetilde{b}_n\}_1^N$ are constants from $\mathbb{C}$ (Remark \ref{r3.9}). Using (ii)\eqref{eq3.16} relative to $s_2(\lambda)$, we construct $s_3(\lambda)$ ($=s_2^{-1}(\lambda\zeta_3)$).

Knowing $s_2(\lambda)$, we construct function $\chi(\lambda,x)$ \eqref{eq3.43}, and from $s_2(\lambda)$, $s_3(\lambda)$ we find the functions $G_2(\lambda,x)$ and $G_3(\lambda,x)$ \eqref{eq3.48}. Finally, define the function $b(\lambda,x)$ where $\{\varkappa_n\}_1^N$ are from \eqref{eq4.1}, and $R_n(x)$ and $\widetilde{R}_n(x)$ are unknown. In terms of the functions
\begin{equation}
\varphi_2(t,x)\stackrel{\rm def}{=}c_3(-i\zeta_3t,x)\psi_2(-it,x);\quad\varphi_3(t,x)\stackrel{\rm def}{=}c_2(-i\zeta_2t,x)\psi_3(-it,x),\label{eq4.2}
\end{equation}
system \eqref{eq3.50} becomes
\begin{equation}
\left\{
\begin{array}{cccc}
{\displaystyle\varphi_2(t,x)Q_3(t,x)+\frac{\zeta_3}{2\pi i}\int\limits_0^\infty\frac{d\tau}{\tau-\zeta_3t}\varphi_2(\tau,x)-\frac{\zeta_2}2\varphi_3(t,x)}\\
{\displaystyle-\frac{\zeta_2}{2\pi i}/\hspace{-4.4mm}\int\limits_0^\infty\frac{d\tau}{\tau-t}\varphi_3(\tau,x)-b(-i\zeta_2t,x)=1;}\\
{\displaystyle\varphi_3(t,x)Q_2(t,x)-\frac{\zeta_2}{2\pi i}\int\limits_0^\infty\frac{d\tau}{\tau-\zeta_2t}\varphi_3(\tau,x)+\frac{\zeta_3}2\varphi_2(t,x)}\\
{\displaystyle+\frac{\zeta_3}{2\pi i}\int\limits_0^\infty\frac{d\tau}{\tau-t}\varphi_2(\tau,x)-b(-i\zeta_3t,x)=1}
\end{array}\right.\label{eq4.3}
\end{equation}
where
\begin{equation}
Q_1(t,x)\stackrel{\rm def}{=}c_2^{-1}(-i\zeta_2t,x)\chi^{-1}(-i\zeta_3t,x);\quad Q_3(t,x)\stackrel{\rm def}{=}c_3^{-1}(-i\zeta_3t,x)\chi^{-1}(-i\zeta_2t,x).\label{eq4.4}
\end{equation}
Equalities \eqref{eq3.51}, \eqref{eq3.52}, in view of \eqref{eq4.2}, become
\begin{equation}
\left\{
\begin{array}{cccc}
{\displaystyle\frac{\zeta_3}{2\pi i}\int\limits_0^\infty\frac{d\tau}{\tau-i\zeta_3\varkappa_m}\varphi_2(\tau,x)-\frac{\zeta_2}{2\pi i}\int\limits_0^\infty\frac{d\tau}{\tau-i\varkappa_m}\varphi_3(\tau,x)-\theta_m(x)R_m(x)}\\
-b(\zeta_2\varkappa_m,x)=1\quad(1\leq m\leq N);\\
{\displaystyle\frac{\zeta_3}{2\pi i}\int\limits_0^\infty\frac{d\tau}{\tau+i\varkappa_m}\varphi_2(\tau,x)-\frac{\zeta_2}{2\pi i}\int\limits_0^\infty\frac{d\tau}{\tau+i\zeta_2\varkappa_m}\varphi_3(\tau,x)-\widetilde{\theta}_m(x)\widetilde{R}_m(x)c}\\
-b(-\zeta_3\varkappa_m,x)=1\quad(1\leq m\leq N)
\end{array}\right.\label{eq4.5}
\end{equation}
where
\begin{equation}
\theta_m(x)\stackrel{\rm def}{=}\zeta_3e^{-\sqrt3\zeta_3\varkappa_mx}b_m\frac{\chi(\varkappa_m,x)}{\chi(\zeta_2\varkappa_m,x)};\quad\widetilde{\theta}_m(x)\stackrel{\rm def}{=}\zeta_2e^{-\sqrt3\zeta_2\varkappa_mx}\widetilde{b}_m\frac{\chi(-\varkappa_m,x)}{\chi(-\zeta_3\varkappa_m,x)}\label{eq4.6}
\end{equation}
($1\leq m\leq N$, and the numbers $b_m$, $\widetilde{b}_m$ are from the set \eqref{eq4.1} (Remark \ref{r3.9})). Upon substituting solution $\varphi_2(t,x)$, $\varphi_3(t,x)$, $\{R_n(x)\}_1^N$, $\{\widetilde{R}_n(x)\}_1^N$ to the system \eqref{eq4.3}, \eqref{eq4.5} into formula \eqref{eq3.49}, we obtain a holomorphic in the sector $\Omega_2$ function $\psi_1(\lambda,x)$. Afterwards, using (a) \eqref{eq2.34}, we find
\begin{equation}
\lim\limits_{x\rightarrow\infty}3\lambda^2(\psi_1(\lambda,x)-1)=i\int\limits_x^\infty q(t)dt\quad(\lambda\in\Omega_2),\label{eq4.7}
\end{equation}
hence, after differentiation, we obtain the potential $q(x)$.

So, solution to the inverse problem is equivalent to solvability of the linear system of $2(N+1)$ equations relative to $\varphi_2(t,x)$, $\varphi_3(t,x)$, $\{R_n(x)\}_1^N$, $\{\widetilde{R}_n(x)\}_1^N$.
\vspace{5mm}

{\bf 4.2} Let
$$\vec{R}(x)=\col[R_1(x),...,R_N(x),\widetilde{R}_1(x),...,\widetilde{R}_N(x)];\quad\vec{e}\stackrel{\rm def}{=}\col[1,...,1],$$
then equation system \eqref{eq4.5} in matrix form becomes
\begin{equation}
(A(x)+D(x))\vec{R}(x)=\vec{\varphi}_2(x)-\vec{\varphi}_3(x)-\vec{e},\label{eq4.8}
\end{equation}
here matrix functions $A(x)$ and $D(x)$ are
\begin{equation}
\begin{array}{ccc}
{\displaystyle A(x)\stackrel{\rm def}{=}\left[
\begin{array}{ccc}
z_k(\zeta_2\varkappa_s,x)&\widetilde{z}_k(\zeta_2\varkappa_s,x)\\
z_k(-\zeta_3\varkappa_s,x)&\widetilde{z}_k(-\zeta_3\varkappa_s,x)
\end{array}\right];}\\
D(x)\stackrel{\rm def}{=}\diag[\theta_1(x),...,\theta_N(x),\widetilde{\theta}_1(x),...,\widetilde{\theta}_N(x)],
\end{array}\label{eq4.9}
\end{equation}
and the vector-valued functions $\vec{\varphi}_2(x)$, $\vec{\varphi}_3(x)$ are given by the formulas
\begin{equation}
\begin{array}{cccc}
{\displaystyle\vec{\varphi}_2(x)\stackrel{\rm def}{=}\frac{\zeta_3}{2\pi i}\int\limits_0^\infty d\tau\varphi_2(\tau,x)\col[(\tau-i\zeta_3\varkappa_1)^{-1},...,(\tau-i\zeta_3\varkappa_N)^{-1},(\tau+i\varkappa_1)^{-1},...,}\\
(\tau+i\varkappa_N)^{-1},];\\
{\displaystyle\vec{\varphi}_3(x)\stackrel{\rm def}{=}\frac{\zeta_2}{2\pi i}\int\limits_0^\infty d\tau\varphi_3(\tau,x)\col[(\tau-i\varkappa_1)^{-1},...,(\tau-i\varkappa_N)^{-1},(\tau+i\zeta_2\varkappa_1)^{-1},...,}\\
(\tau+i\zeta_2\varkappa_N)^{-1}].
\end{array}\label{eq4.10}
\end{equation}
If $b_n\not=0$ and $\widetilde{b}_n\not=0$ ($\forall n$), then the matrix-valued function $A(x)+D(x)$ is invertible for $x\gg1$ since the modules of $\theta_n(x0$ and $\widetilde{\theta}_n(x)$ are larger than the modules of matrix elements of $A(x)$. Relation \eqref{eq4.8} implies that
$$\vec{R}(x)=(A(x)+D(x))^{-1}(\vec{\varphi}_2(x)-\vec{\varphi}_3(x)-\vec{e});$$
and due to \eqref{eq3.56}
$$b(\lambda,x)=\langle(A(x)+D(x))^{-1}(\vec{\varphi}_2(x)-\vec{\varphi}_3(x)-\vec{e}),\overline{\vec{z}}(\lambda,x)\rangle_{E^{2N}}$$
where the inner product is taken in $E^{2N}$ and $\overline{\vec{z}}(\lambda,x)$ is a column vector derived from
$$\vec{z}(\lambda,x)=\col[z_1(\lambda,x),...,z_N(\lambda,x),\widetilde{z}_1(\lambda,x),...,\widetilde{z}_N(\lambda,x)]$$
upon complex conjugation of its elements ($z_k(\lambda,x)$, $\widetilde{z}_k(\lambda,x)$ are given by \eqref{eq3.54}). Thus, the function
$$b(\lambda,x)=\langle(A(x)+D(x))^{-1}\vec{\varphi}_2(x),\overline{\vec{z}(\lambda,x)}\rangle-\langle(A(x)+D(x))^{-1}\vec{\varphi}_3(x),\overline{{\vec{z}(\lambda,x)}}\rangle$$
$$-\langle(A(x)+D(x))^{-1}\vec{e},\overline{\vec{z}(\lambda,x)}\rangle,$$
up to the summand $\langle(A(x)+D(x))^{-1}e,\overline{\vec{z}(\lambda,x)}\rangle=b_0(\lambda,x)$, is a linear functional of $\varphi_2(t,x)$ and $\varphi_3(t,x)$  due to \eqref{eq4.10}, therefore
\begin{equation}
b(\lambda,x)=\int\limits_0^\infty\varphi_2(\tau,x)K_2(\tau,x,\lambda)d\tau+\int\limits_0^\infty\varphi_3(\tau,x)K_3(\tau,x,\lambda)d\tau+b_0(\lambda,x).\label{eq4.11}
\end{equation}
Besides, the functions $K_2(\tau,x,\lambda)$ and $K_3(\tau,x,\lambda)$, for every $x\in\mathbb{R}_+$ and $\lambda\in\mathbb{C}$, belong to $L^2(\mathbb{R}_+)$ (relative to $\tau$). Using representation \eqref{eq4.11} of the function $b(\lambda,x)$, re-write system \eqref{eq4.3} in the matrix form:
\begin{equation}
\begin{array}{ccc}
{\displaystyle\left[
\begin{array}{ccc}
 Q_3(t,x)&{\displaystyle-\frac{\zeta_2}2}\\
{\displaystyle\frac{\zeta_3}2}&Q_2(t,x)
\end{array}\right]\vec{\phi}(t,x)+\left[
\begin{array}{ccc}
0&-\zeta_2\\
\zeta_3&0
\end{array}\right]\frac1{2\pi i}/\hspace{-4.4mm}\int\limits_0^\infty\frac{d\tau}{\tau-t}\vec{\phi}(\tau,x)d\tau}\\
{\displaystyle+\int\limits_0^\infty N(\tau,t,x)\vec{\phi}(\tau,x)d\tau=\vec{\phi}_0(t,x),}
\end{array}\label{eq4.12}
\end{equation}
where
$$\vec{\phi}(t,x)\stackrel{\rm def}{=}\col[\varphi_2(t,x),\varphi_3(t,x)];\quad\vec{\phi}_0(t,x)\stackrel{\rm def}{=}\col[b_0(-i\zeta_2t,x)+1,b_0(-i\zeta_3t,x)+1],$$
and the matrix-valued function $N(\tau,t,x)$ is
\begin{equation}
N(\tau,t,x)\stackrel{\rm def}{=}\left[
\begin{array}{ccc}
{\displaystyle\frac{\zeta_3}{2\pi i(\tau-\zeta_3t)}-K_2(\tau,x,-i\zeta_2t)}&-K_3(\tau,x,-i\zeta_2t)\\
-K_2(\tau,x,-i\zeta_3t)&{\displaystyle-\frac{\zeta_2}{2\pi i(\tau-\zeta_2t)}-K_3(\tau,x,-i\zeta_3t)}
\end{array}\right].\label{eq4.13}
\end{equation}
Upon multiplying equality \eqref{eq4.13} from the left by the matrix
$$\gamma\stackrel{\rm def}{=}\left[
\begin{array}{ccc}
0&\zeta_2\\
-\zeta_3&0
\end{array}\right],$$
we obtain
\begin{equation}
P(t,x)\vec{\phi}(t,x)+\frac1{2\pi i}/\hspace{-4.4mm}\int\limits_0^\infty\frac{d\tau}{\tau-t}\vec{\phi}(\tau,x)+\int\limits_0^\infty\widehat{N}(\tau,t,x)\vec{\phi}(\tau,x)d\tau=\vec{\widehat{\phi}}_0(t,x),
\label{eq4.14}
\end{equation}
where
\begin{equation}
P(t,x)\stackrel{\rm def}{=}\left[
\begin{array}{ccc}
{\displaystyle\frac12}&\zeta_2Q_2(t,x)\\
-\zeta_3Q_3(t,x)&{\displaystyle\frac12}
\end{array}\right];\quad\vec{\widehat{\phi}}_0(t,x)\stackrel{\rm def}{=}\gamma\vec{\phi}_0(t,x).\label{eq4.15}
\end{equation}
Singular matrix equations \eqref{eq4.14} are solved using Riemann boundary value problem \cite{25}. So, assuming
\begin{equation}
H(z,x)\stackrel{\rm def}{=}\frac1{2\pi i}\int\limits_0^\infty\frac{d\tau}{\tau-z}\vec{\phi}(\tau,x),\label{eq4.16}
\end{equation}
relation \eqref{eq4.14} yields that
$$(P(t,x)+I)H^+(t,x)=(P(t,x)-I)H^-(t,x)-\int\limits_0^\infty\widehat{N}(\tau,t,x)\vec{\phi}(\tau,x)d\tau+\vec{\widehat{\phi}}_0(t,x).$$
Taking into account that ${\displaystyle\det(P(t,x)\pm I)=\left(1\pm\frac12\right)^2+1\not=0}$, because $Q_2(t,x)$ $\times Q_3(t,x)\equiv1$, we find that
\begin{equation}
H^+(t,x)=B(t,x)H^-(t,x)-\int\limits_0^\infty\widehat{N}_1(\tau,t,x)\vec{\phi}(\tau,x)+\vec{\phi}_1(t,x),\label{eq4.17}
\end{equation}
here
$$B(t,x)\stackrel{\rm def}{=}(P(t,x)+I)^{-1}(P(t,x)-I)=\frac1{13}\left[
\begin{array}{ccc}
1&8\zeta_2Q_2(t,x)\\
-8\zeta_3Q_3(t,x)&1
\end{array}\right];$$
$$\widehat{N}_1(t,x)\stackrel{\rm def}{=}(P(t,x)+I)^{-1}\widehat{N}(t,x);\quad\vec{\phi}_1(t,x)\stackrel{\rm def}{=}(P(t,x)+I)^{-1}\vec{\widehat{\phi}}_0(t,x);$$
and $H^{\pm}(t,x)=H(t\pm i0,x)$ are boundary values of integral \eqref{eq4.16}. Solution of the Riemann boundary value problem \cite{25} gives new singular integral equation of form \eqref{eq4.14} in which the third summand in the left-hand side is absent. Again, solving the obtained equation, we find $\vec{\phi}(t,x)$. So, solvability of equation systems \eqref{eq4.3}, \eqref{eq4.5} is reduced to the invertibility of the matrix-valued function $(A(x)+D(x))$ ($A(x)$, $D(x)$ are given by \eqref{eq4.9}) and to the solution of singular matrix integral equation \eqref{eq4.14}.
\vspace{5mm}

{\bf 4.3} As {\bf scattering data} of inverse problem, we take the totality
\begin{equation}
(s_2(x),c_1(\lambda),\{\varkappa_n\}_1^N,\{b_n\}_1^N,\{\widetilde{b}_n\}_1^N),\label{eq4.18}
\end{equation}
description of which is as follows.

(i) The set $\{\varkappa_n\}_1^N$ is formed by arbitrary different positive numbers, $0<\varkappa_1<...<\varkappa_N<\infty$ ($0<N<\infty$).

(ii) The totalities $\{b_n\}_1^N$, $\{\widetilde{b}_n\}_1^N$ are composed of any non-zero constants from $\mathbb{C}$ ($b_n\not=0$, $\widetilde{b}_n\not=0$, $\forall n$).

(iii) The function $s_2(\lambda)$ (scattering coefficient) is given on the bundle of straight lines $L$ \eqref{eq1.20} and is holomorphic inside the circle $\mathbb{D}_r=\{\lambda\in\mathbb{C}:|\lambda|<r\}$ ($0<r<\infty$), besides, (i) \eqref{eq3.16} takes place.

(iv) The function $c_1(\lambda)$ (matching coefficient) is holomorphic in $\Omega_2^-\cup\mathbb{D}_r$ and has simple poles at the points $\zeta_3\varkappa_n$, $-\zeta_2\varkappa_n$ from $\Omega_2^-$ ($1\leq n\leq N$). (Cf. \eqref{eq3.11} where $c_1(\lambda)=3\lambda^2\alpha\overline{\theta}/(\alpha e_1^{*'}(\lambda,0)-\beta e_1^*(\lambda,0))$).

(v) For $s_2(\lambda)$ and $c_1(\lambda)$, the following identity holds (unitarity condition \eqref{eq3.23} for the scattering problem):
\begin{equation}
\zeta_3s_2(\lambda)s_3^*(\lambda)+\zeta_2s_3(\lambda)s_2^*(\lambda)+1=\frac{c_1(\lambda)c_1^*(\lambda)}{3\lambda^2}\quad(\forall\lambda\in\mathbb{D}_r)\label{eq4.19}
\end{equation}
where $s_3(\lambda)=s_2^{-1}(\lambda\zeta_3)$ (due to (ii) \eqref{eq3.16}).

\begin{theorem}\label{t4.1}
Let totality \eqref{eq4.18}, for which conditions (i) -- (v) hold, be given. Then there exists such operator $\mathcal{L}_q$ \eqref{eq2.2} -- \eqref{eq2.4}, for which set \eqref{eq4.18} is scattering data.
\end{theorem}
\renewcommand{\refname}{ \rm 
\end{Large}
\end{document}